\input amstex
\documentstyle{specgeom}

\def\tslint{{\int\!\!\!\!\!@!@!@!{-}}}
\def\slint{{\int\!\!\!\!\!\!@!@!@!@!{-}}}
\def\d{d\!{}^{\!\text{\rm--}}\!}

\def\crm{\overline{\Bbb R}_-}

\def\ang#1{\langle {#1} \rangle}

\def\Zfrac{\tsize\frac1{\raise 1pt\hbox{$\scriptstyle z$}}}
\def\zfrac{\frac1{\raise 1pt\hbox{$\scriptscriptstyle z$}}}

\def\rmi{ \Bbb R_-}
\define\tr{\operatorname{tr}}

\define\Tr{\operatorname{Tr}}
\define\TR{\operatorname{TR}}

\document

\topmatter
\title
{A Resolvent Approach to Traces and Zeta Laurent Expansions}
\endtitle
\author{Gerd Grubb}
\endauthor
\address
Copenhagen Univ. Math. Dept.,
Universitetsparken 5, DK-2100 Copenhagen, Denmark.
\endaddress
\email  grubb\@math.ku.dk \endemail
\rightheadtext{Resolvent approach}
\thanks
{This work was supported in part by The Danish Science Research Council,
SNF grant 21-02-0446}
\endthanks

\abstract
Classical pseudodifferential operators $A$ on closed
manifolds are considered. It is shown that the basic properties of
the canonical 
trace $\TR A$ introduced by Kontsevich and Vishik are easily proved by
identifying it with the leading nonlocal coefficient $C_0(A,P)$ in
the trace expansion of $A(P-\lambda )^{-N}$ (with an auxiliary elliptic
operator $P$), as determined in a joint work with Seeley 1995. The
definition of $\TR A$ is extended from the cases of noninteger order,
or integer order and even-even parity on odd-dimensional manifolds,
to the case of even-odd parity on even-dimensional manifolds.

For the generalized zeta function
$\zeta (A,P,s)=\Tr(AP^{-s})$, extended meromorphically to $\Bbb C$,
$C_0(A,P)$ equals the coefficient of $s^0$ in the Laurent expansion
at $s=0$ when $P$ is invertible. In the mentioned parity cases,
$\zeta (A,P,s)$ is regular at all integer points. The higher Laurent
coefficients $C_j(A,P)$ at $s=0$ are described as leading nonlocal
coefficients $C_0(B,P)$ in 
trace expansions of resolvent expressions $B(P-\lambda )^{-N}$, with $B$
log-polyhomogeneous as defined by Lesch 
(here $-C_1(I,P)=C_0(\log P, P)$ gives the
zeta-determinant). $C_0(B,P)$ is shown to be a quasi-trace in
general, a canonical trace $\TR B$ in restricted cases, and the
formula of Lesch for $\TR B$ in terms of a 
finite part integral of the symbol is extended to the parity cases.

{\it The paper has appeared in AMS Contemp.\
Math.\ Proc.\ 366 (2005), 43--64. The present version includes minor
corrections.}  
\endabstract
 \endsubjclass

\endtopmatter

\head Introduction \endhead
\medskip
The {\it noncommutative residue} $\operatorname{res} A$, the {\it canonical
trace} $\operatorname{TR}A$  and the
{\it zeta-re\-gu\-la\-ri\-zed determinant} 
$\log\det A$ are three constants associated
with the classical 
pseudodifferential operators ($\psi $do's) $A$ on an $n$-dimensional
closed manifold $X$, under 
various hypotheses (Wodzicki  \cite{W},
Guillemin  \cite{Gu}, Kontsevich and Vishik \cite{KV}, Ray and Singer
\cite{RS}). When $P$ is an invertible  
elliptic classical 
$\psi $do on $X$ of order $m>0$ and with spectrum in a subsector
of $\Bbb C$, one can define the 
generalized zeta function $\zeta (A,P,s)$ as the meromorphic
extension of $\Tr(AP^{-s})$ (defined for large $\operatorname{Re}s$)
to $\Bbb C$; it has a simple pole at $s=0$:$$
\zeta (A,P,s)\sim\frac 1s C_{-1}(A,P)+C_0(A,P)+\sum_{j\ge 1}C_j(A,P)s^j.
\tag 0.1$$
Then $$\aligned
C_{-1}(A,P)&=\tfrac1m \cdot
\operatorname{res}A,\\
C_0(A,P)&=\operatorname{TR}A\;\text{
in restricted cases},\\
 C_1(I,P)&=-\log \det P.
\endaligned\tag 0.2$$

We shall investigate these constants, in particular $C_0(A,P)$ and
$C_1(A,P)$,  by use of the knowledge of the structure of resolvent
expressions $A(P-\lambda )^{-N}$. It 
turns out that it is rather easy to show the trace properties of
$C_0(A,P)$ using Grubb and Seeley \cite{GS1, Th\. 2.1 and 2.7},
for the known cases where $A$ is of low order, noninteger order or
``odd class'' 
(with $n$ odd), as well as
for a new case with another parity property. (Section 1.)

$\operatorname{TR}A$ does not extend to  general 
operators of integer order $\nu $, but here
$C_0(A,P)$ can be viewed as a {\it quasi-trace}, in the sense that it
is determined from $A$ modulo local contributions (from the first
$\nu +n+1$ homogeneous symbol terms in $A$ and $P$), and vanishes on
commutators modulo local contributions. The value of $C_0(A,P)$
modulo local terms is a finite part integral of the symbol of $A$, in
local coordinates. (Section 2.)
(The deviation of $C_0(A,P)$ from being an independently defined trace
of $A$ is further studied in  \cite{KV}, Okikiolu \cite{O1},
Melrose and Nistor \cite{MN}; the latter call $C_0(A,P)$ a regularized trace.
See also Cardona, Ducourtioux, Magnot, Paycha \cite{CDMP},
\cite{CDP}, where it is called a weighted trace.)

In the study of $C_1(A,P)$ and the $C_j(A,P)$ with higher $j$ one
meets the necessity of considering resolvent expressions where $A$ is
replaced by a  log-polyhomogeneous $\psi $do $B$. Here we can use the
results of Lesch \cite{L} to extend the canonical trace
$\operatorname{TR}$ to such operators $B$. This was done for cases
with noninteger or low order in \cite{L}; we now include also higher
integer order cases with 
parity properties, and show that $C_0(B,P)$ is in general a
quasi-trace. In particular, we can identify $\log \det P$ and 
higher derivatives of $\zeta (I,P,s)$ at $s=0$ as quasi-traces;
canonical traces in particular situations. (Section 3.)

Our method relies on an analysis of integrals of symbols, and 
involves neither comparison of meromorphic extensions nor homogeneous
distributions. It moreover allows
us to extend the explicit formula of Lesch \cite{L} for the density
$\omega _{\operatorname{TR}}(B)$ defining 
$\operatorname{TR}B$, to the new integer-order cases. 
--- The strategy is useful in situations where complex powers of
operators are not easy to study directly; for example in
our treatment of boundary value problems jointly with Schrohe \cite{GSc}.

\bigskip
\head 1. The canonical trace  \endhead 
\medskip
Consider an $n$-dimensional compact $C^\infty $
manifold $X$ without boundary. We denote $\{0,1,2,\dots\}=\Bbb N$.
Let $A$ be a classical (i.e., one-step polyhomogeneous) $\psi $do of
order $\nu \in\Bbb R$, acting on the sections of a $C^\infty $ vector
bundle $E$ over $X$. Let $P$ be a
classical elliptic $\psi $do of positive integer order $m$,
likewise acting in $E$ and such that the principal symbol has no
eigenvalues on $\Bbb R_-$. It is shown
in \cite{GS1, Th\. 2.7} by use of calculations in local coordinates
that the operator family
$A(P-\lambda )^{-N}$ (for $N>(n+\nu )/m$) has 
an asymptotic expansion of the trace:
$$
\Tr\bigl( A(P-\lambda )^{-N}\bigr)
\sim
\sum_{ j\in \Bbb N  } \tilde c_{ j}(-\lambda ) ^{\frac{\nu +n -j}m-N}+ 
\sum_{k\in \Bbb N}\bigl( \tilde c'_{ k}\log (-\lambda ) +\tilde c''_{
k}\bigr)(-\lambda ) ^{{ -k}-N},
\tag1.1$$
for $\lambda \to\infty $ on rays in 
an open subsector of $\Bbb
C$ containing $\Bbb R_-$.
\comment
Here each $\tilde c_j$ (and each $\tilde c_k'$) comes from a specific
homogeneous term in the symbol of $A(P-\lambda)^{-N}$, whereas the 
$\tilde c''_k$ depend on the full symbol. So the
coefficients $\tilde c_j$ and $\tilde c_k'$ depend each on a finite
set of homogeneous terms in the symbols of $A$ and $P$; we call such
coefficients `locally determined' (or `local'), while the $\tilde
c''_k$ are called `global'. When $\nu \notin
\Bbb Z$, the $\tilde c'_k$ vanish. 

\endcomment 

In  local coordinates, the term of degree $\nu -mN-j$ in the
symbol of \linebreak$A(P+\mu ^m)^{-N}$ determines the coefficient $\tilde c_j$
and, if $\nu \in\Bbb Z$ and $\frac{j- \nu -n}m\in\Bbb N$, the
coefficient $\tilde c'_k$ with 
$k=\frac{j-\nu -n}m$ (one sets the $\tilde c'_k$ with $mk+\nu +n<0$ equal to
0). If $\nu \notin\Bbb Z$, $\tilde c'_k=0$ for all $k\ge
0$.
In terms of the original operators, $\tilde c_j$ and (if $\nu \in\Bbb
Z$ and $\frac{j-\nu -n}m\in\Bbb N$) $\tilde c'_{(j-\nu -n)/m}$ depend
solely on the (strictly homogeneous part of the)
homogeneous terms of 
degrees $\nu ,\nu -1,\dots,\nu -j$ resp\. $m,m-1,\dots,m-j$ in the
symbols of $A$ resp\. $P$ (in short: the first $j+1$ homogeneous
terms). As in \cite{GS1} and many other works, we call such coefficients
``locally determined'' or just ``local''. The coefficients $\tilde
c''_k$ depend on  
the full structure of the operators on the manifold (are ``global'').

Note that when $\nu \in\Bbb Z$ and $\frac{j-n-\nu}m$ is an integer
$k\ge 0$, both $\tilde c_j$ and $\tilde c''_k$ 
contribute to the power $(-\lambda )^{-k-N}$. Their sum is
independent of the choice of local coordinates, whereas the
splitting in $\tilde c_j$ and $\tilde c''_k$ depends in a well-defined way
on the symbol structure in the chosen local coordinates (see  \cite{GS1,
Th.\ 2.1} or the elaboration in Theorem 1.3 below).    

Along with (1.1) there is the essentially equivalent expansion
(the transition between (1.1) and (1.2) is
accounted for e.g\. in \cite{GS2}):
$$
\Gamma (s)\Tr (AP^{-s})\sim
\sum_{ j\in \Bbb N } \frac{c_{j}}{s+\frac{j-\nu -n}m}-\frac{\Tr (A\Pi
_0(P))}s +  \sum_{k\in \Bbb N}\Bigl(\frac{
c'_{k}}{(s+k)^2}+\frac{c''_{k}}{s+ k}\Bigr) .\tag1.2
$$ 
This means that $\Gamma (s)\Tr (AP^{-s})$,
defined in a standard way for $\operatorname{Re}s>\frac {\nu +n}m$, 
extends meromorphically to $\Bbb C$ with the pole structure 
indicated in the right hand side. Here $\Pi _0$ is the orthogonal
projection onto
the nullspace of $P$ (on which $P^{-s}$ is taken to be zero).
The coefficients $\tilde c_j$ and
$c_j$, resp.\ $\tilde c'_k$ and $ c'_k$,
are proportional by universal nonzero
constants. When the $c'_k$ vanish (e.g., when $\nu +n\notin \Bbb N$),
the same holds for $\tilde c''_k$ and $ c''_k$. More generally,
the pair $\{\tilde c'_k, \tilde c''_k\}$ is for each $k$
universally related to the pair $\{c'_k, c''_k\}$ in a linear way.
In particular, $\tilde c'_0=c'_0$, and $\tilde
c''_0=c''_0$ if $c'_0=0$, and when $\nu $ is an integer $\ge -n$,
$\tilde c_{\nu 
+n}= c_{\nu +n}$.
We shall {\it define}
$$\tilde c_{\nu +n}= c_{\nu +n}=0\text{ if }\nu <-n\text{
or }\nu \notin \Bbb Z;\tag1.3$$
 then the identifications hold in these cases too. 

We are particularly interested in
$C_0(A,P)$, defined by
$$
C_0(A,P)=c_{\nu +n}+ c''_0, \text{ equal to }\tilde c_{\nu +n}+\tilde
c''_0\text{ if }N=1.
\tag1.4
$$
When $N=1$, $C_0(A,P)$ is the coefficient of $(-\lambda )^{-1}$ in (1.1).
For general $N$, the coefficient of $(-\lambda )^{-N}$ in (1.1) satisfies
$ \tilde c_{\nu +n}+\tilde
c''_0=C_0(A,P)-\alpha _Nc'_0$, 
where $\alpha _N=\sum _{1\le j<N}\frac1{j}$, cf.\ \cite{G5, Lemma 2.1}.
[The preceding lines, modified in November 2005, correct the printed version
of the present 
paper, where the term with $\alpha _N$ was missing. However, in the
sequel, $C_0(A,P)$ is determined in cases where $c'_0=0$, so the
results in the following remain valid.]

Division by $\Gamma (s)$ in (1.2) gives
the structure of the meromorphic extension of $\Tr(AP^{-s})$,
also denoted $\zeta (A,P,s)$. When $\Pi _0(P)=0$, it
has the Laurent expansion at $s=0$:
$$
\zeta (A,P,s)\sim \frac 1s {C_{-1}(A,P)}+C_0(A,P)+\sum _{l\ge
1}C_{l}(A,P)s^l,\text{ with }C_{-1}(A,P)= {c'_0};\tag1.5
$$
here $C_0(A,P)$ must be replaced by $C_0(A,P)-\Tr(A\Pi _0(P))$ if
$\Pi _0(P)\ne 0$.

If the eigenvalues of the principal symbol of $P$ lie in a sector
$\{\lambda \mid |\arg\lambda|\le \theta \}$ with $0\le \theta
<\frac\pi 2$, so that $e^{-tP}$ is well-defined, there is a third
trace expansion that is equivalent with (1.1) and (1.2) (cf\. e.g\.
\cite{GS2} for the transition):
$$
\Tr  (Ae^{-tP})\sim
\sum_{ j\in \Bbb N } c_{j}t^{\frac{j-\nu -n}m}+ 
\sum_{k\in \Bbb N}\bigl({- c'_{k}}\log t+{c''_{k}}\bigr) t^{{k}},\tag1.6
$$
for $t\to 0+$; the coefficients here are {\it the same} as those in (1.2).

One interest of studying the resolvent-type expansion
(1.1) along with
(1.2) is that it allows to determine the coefficients from specific
integrals of symbols (cf.\ \cite{GS1, Th\. 2.1}, Th\. 1.3 below).  
In the following, we use the notions of \cite{GS1} without  taking space
up with repetition of basic rules
of calculus explained there.

In \cite{W}, Wodzicki introduced a trace functional on the full algebra of
classical $\psi $do's, vanishing on trace-class operators; it is
usually denoted $\operatorname{res}A$ and is called {\it the
noncommutative residue of} $A$. In the above situation, it satisfies$$
\operatorname{res} A=m\cdot \operatorname{Res}_{s=0}\Tr (AP^{-s})
=m\cdot  c'_0=m\cdot \tilde c'_0.\tag1.7
$$
See also Guillemin \cite{Gu} and
the survey of Kassel \cite{K}.

In \cite{KV}, Kontsevich and Vishik introduced a different trace functional 
$\operatorname{TR}A$, called {\it the canonical trace}, which extends
the standard trace for 
trace-class operators, but is only
defined for part of the higher-order $\psi $do's.
We shall show that the following definition is
consistent with that of
\cite{KV}:

\proclaim{Definition 1.1} Let $A$ be a classical $\psi $do in $E$ of
order $\nu \in\Bbb R$, and let $P$ be a classical elliptic $\psi $do
in $E$ of even order $m>0$ such that the principal symbol has no eigenvalues
on $\Bbb R_-$. Assume that one of the following statements is
verified (with notation explained around {\rm (1.10)--(1.11)} below):
\roster
\item "(1)" $\nu <-n$.
\item "(2)" $\nu \notin \Bbb Z$.
\item "(3)" $\nu \in\Bbb Z$, 
 $A$ is 
even-even, and $n$ is odd.
\item "(4)" $\nu \in\Bbb Z$, 
$A$ is even-odd, and $n$ is even.
\endroster
In the cases {\rm (3)} and {\rm (4)}, take $P$ even-even.
Then 
$$\operatorname{TR}A=\tilde c''_0 = c''_0=C_0(A,P).\tag1.8
$$ 
\endproclaim 

Definitions for the cases (1), (2) and (3) were given in \cite{KV},
whereas the case (4) is new. The rest of Section 1 will be devoted to the
justification of Definition 1.1. [Added November 2005:] It is seen in
particular that 
$c'_0$ vanishes in all four cases, so that $C_0(A,P)$ identifies
directly with the coeficient of  $(-\lambda )^{-N}$ in (1.1).

In case (1), the definition is consistent with \cite{KV} in view of the
following fact:

\proclaim{Lemma 1.2}
When $\nu <-n$, then $c_{\nu +n}=c'_0=0$ and
$c''_0=\Tr A$.
\endproclaim 

\demo{Proof}
One has in this case that 
$\Tr\bigl(
A(P-\lambda )^{-N}(-\lambda )^N\bigr)=\tilde c''_0+O(\lambda
^{-\varepsilon })$ with $\varepsilon >0$ 
when $\nu <-n$, in view of (1.1) and the comments on the 
coefficients. Then the identity follows since $A(P-\lambda )^{-N}(-\lambda )^N\to A$ in
trace-norm for $\lambda \to -\infty $ on $\Bbb R_-$. 

For a detailed proof of the latter 
fact, observe that $R_\lambda =(P-\lambda )^{-1}$ satisfies$$
R_\lambda ^{N}(-\lambda )^N-1=[(-\lambda )^N-(P-\lambda
)^N](P-\lambda )^{-N}= PR_\lambda M_\lambda ,\tag1.9
$$
where $M_\lambda =(P-\lambda )^{1-N}\sum_{0\le j\le N-1}\binom N j(-\lambda
)^jP^{N-1-j}$ is uniformly bounded in $L^2(X,E)$ operator norm for $\lambda 
\le -1$. Then with $\delta =\min\{\frac1{2m}(-n-\nu ),1\}>0$, $$
\multline
\|(AR_\lambda ^{N}(-\lambda )^N- A)f\|_{H^{n+m\delta }} \le 
c\|(R_\lambda ^{N}(-\lambda )^N- I)f\|_{H^{-m\delta }} 
\\
\le 
c'|\lambda
|^{-\delta } |\lambda |^{\delta }\|R_\lambda M_\lambda
f\|_{H^{m-m\delta }} \le c'' |\lambda |^{-\delta 
}\|f\|_{L^2},
\endmultline
$$
using that $|\lambda |^{\delta }\|R_\lambda f\|_{H^{m-m\delta
}}\le c_3(\|R_\lambda f\|_{H^m}+|\lambda |\|R_\lambda f\|_{L^2})\le
c_4\|f\|_{L^2}$. Thus for $\lambda \to -\infty $,$$ 
\|AR_\lambda ^{N}(-\lambda )^N- A\|_{\Tr}\le c_5\|AR_\lambda ^{N}(-\lambda )^N- A\|_{\Cal L(L^2,H^{n+m\delta })}\to 0.\quad\square
$$
\enddemo 

In the case (2), the formula (1.8) makes sense
since $c'_0=0$ and $c_{\nu +n}=0$ by definition.
It was shown by Lesch in \cite{L} that
the definition is consistent with that
of \cite{KV} in this case, if $P$ is selfadjoint positive with scalar
leading symbol. The following analysis shows that the definition is likewise
consistent with that of \cite{KV} for the $P$'s considered here.

We now turn to (3)
and (4): 
As in \cite{G2, Sect\. 5}, we say that a classical $\psi $do $Q$ of order $r\in\Bbb Z$ with symbol $q\sim
\sum_{l\in\Bbb N}q_{r-l}(x,\xi )$ ($q_{r-l}$ $C^\infty $ in $(x,\xi )$
and homogeneous of degree
$r-l$ in $\xi $ for $|\xi |\ge 1$) has {\bf even-even} alternating
parity (in short: is even-even), when the 
symbols with even (resp\. odd) degree $r-l$ are even (resp\. odd) in $\xi $:
$$
q_{r-l}(x,-\xi )=(-1)^{r-l}q_{r-l}(x,\xi ) \text{ for }|\xi |\ge 1.\tag1.10
$$ 
The operator (or symbol) is said to have {\bf even-odd} alternating
parity in the reversed 
situation where the
symbols with even (resp\. odd) degree $r-l$ are odd (resp\. even) in $\xi $:
$$
q_{r-l}(x,-\xi )=(-1)^{r-l-1}q_{r-l}(x,\xi ) \text{ for }|\xi |\ge 1.\tag1.11
$$ 
Other authors use other names, e.g\. \cite{KV} calls the even-even
symbols ``odd-class'', studying them on odd-dimensional manifolds, and
Okikiolu \cite{O3} uses the words ``regular parity'' resp\. ``singular
parity'' for the 
even-even resp\. even-odd alternating parity.
Differential operators and their parametrices are even-even, whereas
e.g.\ 
$|A|=(A^2)^\frac12$ is even-odd, when $A$ is a first-order elliptic
selfadjoint differential operator (as noted in \cite{GS2, p\. 48}).

In case (3), defining $\TR A$ as $C_0(A,P)$ is consistent with
\cite{KV, Sect\. 
7.3}, cf\. also \cite{O3}. We shall now show (from scratch) that
this constant has
the desired properties, and that the constant in case (4) likewise does so.
The proof --- inspired from \cite{G2, Theorem 5.2} --- 
shows that in the cases (3) and (4), the logarithmic terms and the
local terms with $\nu +n-j$ even vanish.
It will be based on an exact application of the method of proof of \cite{GS1,
Th\. 2.1} to the present operator family $A(P-\lambda )^{-N}$. Here
we moreover give an account of how the coefficient
$C_0(A,P)$ looks for
general $A$.

Write $$
A(P-\lambda )^{-N}=A(P+\mu ^m)^{-N}=Q(\mu ), \quad\mu
=(-\lambda )^{\frac1m},\tag 1.12$$ 
where
$\mu $ is included in the symbol as in \cite{G1}, \cite{GS1}, \cite{G2};
then $(P+\mu ^m)^{-N}$ is weakly polyhomogeneous in $(\xi ,\mu )$
(strongly so if $P$ is a differential operator).

There is a finite cover of $X$ by coordinate patches $U_i$ ($i\le
i_0$) with trivializations of $E$; $\psi _i\:E|_{U_i}\to V_i\times
\Bbb C^{\operatorname{dim}E}$ with $V_i\subset\subset \Bbb R^n$,  and a subordinate
partition of 
unity $\varphi _j$ ($j\le j_0$) such that any two of the functions
$\varphi _j$ are supported in one of the $U_i$'s (for
$i=i(j_1,j_2)$).
Then we can write$$
A=\sum_{j_1,j_2\le j_0}\varphi _{j_1}A\varphi _{j_2},\tag 1.13
$$ 
a finite sum of $\psi $do's, each acting in a coordinate patch (and
preserving the property of being supported in the patch). Since
the coefficients in the trace expansions depend linearly on $A$, it
suffices to consider the expansions for each term in (1.13).
Actually, we can, by linear translations in $\Bbb R^n$, replace the
$V_{i(j_1,j_2)}$ by sets $V'_{j_1,j_2}$ with a
positive distance from one another, so that $A$ in (1.13) carries over to a
$\psi $do that is a sum of pieces supported in each $V'_{j_1,j_2}$ for
$j_1,j_2\le j_0$ --- we shall denote it $A$ again.
We likewise
consider $(P-\lambda )^{-N}$ in the coordinate patches
carried over to $\Bbb R^n$ in this way.

In the localized situation, let $Q(\mu )$
have the symbol$$
q(x,\xi ,\mu )\sim \sum_{j\in\Bbb N}q_{\nu -mN-j}(x,\xi ,\mu );\tag1.14
$$
here the $q_{\nu -mN-j}$ are homogeneous
in $(\xi ,\mu 
)$ of degree $\nu -mN-j$, for $|\xi |\ge 1$.
For simplicity of notation, we can let $\mu $ run on the ray $\Bbb
R_+$ (other rays are treated similarly, and holomorphy in $\mu $ is assured by
\cite{GS1, Lemma 2.3}).
Besides this polyhomogeneous structure, the important knowledge is
that the symbol of $Q$ lies in a suitable $S^{k,d}$-space, as defined
in \cite{GS1}. When a symbol $f(x,\xi ,\mu )$ lies in $ S^{k,d}$, the
$d$-index indicates 
that $f(x,\xi ,\mu )=\mu 
^df_1(x,\xi ,\mu )$, where $f_1$ has a Taylor expansion in $z=\frac1\mu
$ at $z=0$:
$$\aligned
f(x,\xi ,\mu )&=\mu ^{d}f_1(x,\xi ,1/z)=z^{-d}\sum_{0\le l<L}f^{(l)}(x,\xi
)z^{l}+O(\ang{\xi 
}^{k+L}z ^{-d+L})\\
&=\sum_{0\le l<L}f^{(l)}(x,\xi )\mu ^{d-l}+O(\ang{\xi
}^{k+L} \mu  ^{d-L}),\text{ for any }L,\endaligned\tag1.15
$$ 
with $f^{(l)}\in S^{k+l}$ (see \cite{GS1} for further details; $\ang \xi
$ stands for $(1+|\xi |^2)^{\frac12}$).

\proclaim{Theorem 1.3}
{\rm (i)} In the localized situation, 
when $N>(\nu +n)/m$, the diagonal value of the kernel $K( Q(\mu
),x,y)$ of $Q(\mu )$ has an
asymptotic expansion$$\multline
K( Q(\mu ),x,x)
\sim
\sum_{ j\in \Bbb N  } \tilde c_{ j}(x)\mu  ^{\nu +n -j-mN}+ 
\sum_{k\in \Bbb N}\bigl(m \tilde c'_{ k}(x)\log \mu  +\tilde c''_{
k}(x)\bigr)\mu ^{-mk-mN}\\
\sim
\sum_{ j\in \Bbb N  } \tilde c_{ j}(x)(-\lambda ) ^{\frac{\nu +n -j}m-N}+ 
\sum_{k\in \Bbb N}\bigl( \tilde c'_{ k}(x)\log (-\lambda )  +\tilde c''_{
k}(x)\bigr)(-\lambda ) ^{-k-N}.
\endmultline\tag1.16
$$
Here, when we define $\tilde c_{\nu +n}(x)=0$ if $\nu <-n$ or $\nu \in \Bbb
R\setminus \Bbb Z$,
$$
\tilde c_{\nu +n}(x)+\tilde c''_0(x)=\slint a(x,\xi )\,\d\xi +\text{ local terms},\tag1.17
$$
where $\tslint a$ is defined from  the symbol $a(x,\xi )\sim
\sum_{j\in\Bbb N}a_{\nu -j}(x,\xi )$ of $A$ by:$$
\multline
\slint a(x,\xi )\d \xi =\sum_{j\le \nu +n}
\bigl(\int_{|\xi |\le 1} a_{\nu -j}(x,\xi )\,\d\xi -\tfrac{1-\delta _{\nu +n,j}}{\nu +n-j}
\int_{|\xi |=1}a_{\nu -j}(x,\xi )\,\d S(\xi )\bigr)\\
+\int_{\Bbb
R^{n}}
\bigl(a(x,\xi )-\sum_{j \le \nu +n}a_{\nu -j}(x,\xi )\bigr)\,\d\xi, 
\endmultline\tag 1.18
$$
 and the local terms depend only on the first $n+[\nu] +1$ strictly
homogeneous 
terms in the symbols of $A$ and $P$. When $\nu \notin \Bbb Z$, all
$\tilde c'_k$ vanish.

{\rm (ii)} If, moreover,  $A$ and $P$ are as in Definition
{\rm 1.1} 
with {\rm (3)} or {\rm (4)}, then the expansion {\rm (1.16)} reduces
to 
the form$$
K(Q(\mu ),x,x)
\sim
\sum_{ j\in \Bbb N ,\,j-n-\nu \text{ odd}} \tilde c_{ j}(x)(-\lambda ) ^{\frac{\nu +n -j}m-N}+ 
\sum_{k\in \Bbb N}\tilde c''_{
k}(x)(-\lambda ) ^{{ -k}-N}.
\tag1.19
$$

{\rm (iii)} In each of the cases {\rm (1)--(4)} in Definition {\rm
1.1}, $\tilde c'_0(x)=0$, $\tilde c_{\nu +n}(x)=0$ (for any choice of
local coordinates), and $$
\tilde c_0''(x)=\slint a(x,\xi )\,\d\xi ,\tag 1.20
$$
clearly independent of $P$.
\endproclaim 

\demo{Proof} In these formulas,
$\delta  _{r,s}$ is the 
Kronecker delta, and we use the notation $[\nu ]$ for the largest integer
$\le \nu $. 
The theorem will be proved by an examination of how the coefficients
in (1.1) arise in the proof of \cite{GS1, Th\. 2.1}. 
In the proof of (ii), we consider for definiteness e.g\. the case (4)
(the other case (3) is treated in a completely analogous
fashion).

It is shown in \cite{GS1, Sect\. 2} that $(P+\mu ^m
)^{-N}$ has symbol in $S^{-mN,0}\cap S^{0,-mN}$, so the symbol 
$q(x,\xi ,\mu )$ of $Q(\mu )=A(P+\mu ^m)^{-N}$ satisfies
$$
q(x,\xi ,\mu )\in S^{\nu -mN,0}\cap S^{\nu ,-mN}.\tag1.21
$$ 
The expansion of $q$ corresponding to (1.15) for the $d$-index equal to
$-mN$ reflects the fact that$$
\aligned
A(P+\mu ^m)^{-N}&=z  ^{mN}A(I+z ^mP)^{-N}
= z
^{mN}A\sum _{0\le l<L}\tbinom {-N}l z ^{ml}P^l+O(z^{m(N+L)})\\
&=\sum_{0\le l<L}\tbinom{-N}l AP^l\mu  ^{-m(N+l)}+O(\mu ^{-m(N+L)}).
\endaligned\tag1.22
$$
In fact, only $m$'th powers enter nontrivially in the expansion of
$q$ (since $q$ is a 
function of $\lambda =-\mu ^m$):
$$
q(x,\xi ,\mu )=\sum_{0\le l<L}q^{(l)}(x,\xi )\mu ^{-m(N+l)}+O(\ang\xi
^{\nu +mL}\mu ^{-m(N+L)}).\tag1.23
$$
Here the $q^{(l)}$ are polyhomogeneous symbols of order $\nu +ml$.
In case (4), $q(x,\xi ,\mu )$ is even-odd and the
$q^{(l)}$ are even-odd.
(It may be observed that the first term in the last sum in (1.22) equals $A\mu
^{-mN}$, so the first coefficient in (1.23) is $q^{(0)}=a(x,\xi )$.)

The kernel of $Q(\mu )=\operatorname{OP}(q)$, restricted to the diagonal
$x=y$, is $$
K_q(x,x,\mu )=\int_{\Bbb R^n}q(x,\xi ,\mu )\,\d\xi,\quad \d\xi =(2\pi
)^{-n}d\xi .\tag1.24 
$$
The contributions from a homogeneous term $q_{\nu -mN-j}$ to (1.16)
are found by 
splitting the corresponding integral in three pieces:
$$
K_{q_{\nu-mN-j}}=\int_{|\xi |\ge \mu }q_{\nu-mN-j}\,\d\xi  +\int_{|\xi |\le 1
}q_{\nu-mN-j}\,\d\xi+\int_{1\le |\xi |\le \mu }q_{\nu-mN-j}\,\d\xi.
\tag1.25 
$$
In the first integral we replace $\xi $ by $\mu \eta $ and use the
homogeneity; this gives a contribution
 to $\tilde c_j(x)\mu ^{\nu +n-j-mN}$.

In case (4), 
$q_{\nu-mN-j}$ is odd in $\xi $ when $\nu -j$ is even, so the
contribution vanishes when $\nu -j$ is even.
Since $n$ is even (as well as $mN$), these are the cases where $\nu
+n-j$ is even.
(In a similar way one sees in case (3) that since $Q$ is
even-even, the 
contributions to cases $\nu +n-j$ even vanish since $n$ is odd.) 

For the other pieces in (1.25) we use moreover that $q_{\nu -mN-j}\in
S^{\nu -j,-mN}$. 
The second piece contributes straightforwardly to the $\tilde c''_k$-terms:
$$
\int_{|\xi |\le 1}q(x,\xi ,\mu )\,\d \xi =\sum_{0\le l<L}\mu
^{-m(N+l)}\int_{|\xi |\le 1}q^{(l)}(x,\xi )\,\d\xi +O(\mu
^{-m(N+L)}),\tag1.26$$
where $q^{(0)}=a$, as noted above.
In the third piece, the terms in the symbol are
homogeneous, and since we are integrating over a 
bounded part of $\Bbb R^n$, we need not worry about integrability at
infinity. We use the expansion (as in (1.23))$$
q_{\nu -mN-j}(x,\xi ,\mu )=\sum_{0\le l<L}q_{\nu -j+ml}^{(l)}(x,\xi )\mu
^{-mN-ml}+R_{j,L}(x,\xi ,\mu ), 
\tag1.27
$$
where the coefficients  $q_{\nu -j+ml}^{(l)}$ are
homogeneous of degree $\nu -j+ml$ in $\xi $ and $R_{j,L}$ is $O(\ang\xi
^{\nu -j+mL} \mu ^{-m(N+L)})$; in case (4), all these terms are even-odd.
One finds by use of polar coordinates:
$$
\aligned
\mu ^{-m(N+l)}&\int_{1\le |\xi |\le \mu }q_{\nu -j+ml}^{(l) }(x,\xi 
)\,\d\xi \\
&=\mu ^{-m(N+l)}\int_1^\mu r^{\nu -j+ml +n-1}\,dr
\int_{|\xi |=1}q_{\nu -j+ml}^{(l)} (x,\xi)\,\d S(\xi )\\
&=\cases \frac{c_{j,l}(x)}{\nu -j+ml+n}(\mu ^{\nu -j-mN +n}-\mu ^{-m(N+l)})&\text{ if
}\nu -j+ml+n\ne 0,\\
c_{j,l}(x)\mu ^{-m(N+l)}\log \mu &\text{ if }\nu -j+ml+n= 0,\endcases
\endaligned\tag1.28
$$
where$$
c_{j,l}(x)= \int_{|\xi |=1}q_{\nu -j+ml}^{(l) }(x,\xi
)\,\d S(\xi ).\tag1.29
$$

When $j\ne \nu
+ml+n$, the term $\frac1{\nu -j+ml+n}c_{j,l}(x)\mu ^{\nu -j-mN+n}$
contributes to the $\tilde c_j(x)$-term, whereas $\frac{-1}{\nu
-j+ml+n}c_{j,l}(x)\mu ^{-m(N+l)}$ is absorbed in the $\tilde
c''_l(x)$-term.
When $j=\nu +ml+n$, we get the $l$'th log-term in the first line of
(1.16) with coefficient
$c_{j,l}(x)$; it will be denoted $m \tilde c'_{l}(x)$ to comply with
the notation conventions of (1.1).
The second line in (1.16) is obtained by insertion of $\mu =(-\lambda
)^{\frac 1m}$. If $\nu \notin\Bbb Z$, logarithmic terms cannot occur.

Let us see how the coefficients look in case (4):
The logarithmic contribution comes when $j=\nu +ml+n$, and then
since $n$ is even and $q_{\nu -j+ml}^{(l)}$ is even-odd of degree
$\nu -j+ml=-n$, $c_{j,l}(x)=0$. 
(Similarly, this coefficient vanishes in case (3) where $q_{\nu
-j+ml}^{(l)}$ is even-even and $n$ is odd.)  
Thus there are no logarithmic terms! Moreover, $c_{j,l}(x)$ 
vanishes if $\nu -j+mN+n$ is even, i.e., when
$\nu
-j+n$ is even, so there is no contribution to $\tilde c_j(x)$ in this
case. Hence the expansion terms from the third piece in
(1.25) only contribute to 
the terms in 
(1.19).

It is accounted for in \cite{GS1, pf\. of Th\. 2.1} (and in more
detail in \cite{GH,
Sect.\ 3}) how the remainders,
from the polyhomogeneous expansion (1.14) as well as the expansions
in powers of $\mu $, are handled; for completeness we recall the
arguments here: 
For the remainder $R_{j,L}$ in (1.27), consider a case where
$\nu -j+mL>0$. As 
noted in \cite{GS1}, $R_{j,L}$ is $O(|\xi |
^{\nu -j+mL} \mu ^{-m(N+L)})$ for $|\xi |\ge 1$, so it extends by
homogeneity for $|\xi |\le 1$ to a continuous function
$R^h_{j,L}(x,\xi ,\mu )$ satisfying the same estimate. Then
$$\aligned
\int_{1\le |\xi |\le \mu }R_{j,L}\,\d\xi &=
\int_{ |\xi |\le \mu }R^h_{j,L}\,\d\xi -
\int_{ |\xi |\le 1 }R^h_{j,L}\,\d\xi \\
&=c'_{j,L}(x)\mu ^{\nu
-j-mN+n}+O(\mu ^{-m(N+L)}),
\endaligned$$ 
giving another contribution to the coefficient $\tilde c_j(x)$. In
the case (4), the
contribution vanishes for $\nu -j+n$ even. 
This shows that the homogeneous terms $q_{\nu -mN-j}$ have expansions
as in (1.16) down to an $O(\mu ^{-m(N+L)})$-error when $L$ is large;
then it holds {\it a fortiori} for small $L$. 
 Now consider the remainder $q'_J=q-\sum_{j<J}q_{\nu -j-mN}$ in the
expansion (1.14); it is in $S^{\nu -J-mN,0}\cap S^{\nu -J,-mN}$
(depending on $\mu $ through $\lambda $), so it has an expansion, for
any $L$,
$$
q'_J(x,\xi ,\mu )=\sum_{0\le l<L}q^{\prime (l)}_J(x,\xi )\mu ^{-m(N+l)}+O(\ang\xi
^{\nu -J+mL}\mu ^{-m(N+L)}),$$
with $q^{\prime (l)}_J(x,\xi )\in S^{\nu -J+ml}$. Assume that $J>\nu
+mL+n$; then all the terms are integrable in $\xi $, and
$$
\int_{\Bbb R^n}q'_J(x,\xi ,\mu )\,\d\xi =\sum_{0\le
l<L}c'_{J,l}(x)\mu ^{-m(N+l)}+O(\mu ^{-m(N+L)}).
$$
The $c'_{J,l}(x)$ are taken into the coefficients $\tilde
c''_l(x)$.
We conclude that there is an asymptotic expansion (1.16) which for
any $L\ge 0$ can be calculated down to an error $O(\mu ^{-m(N+L)})$ by
taking $J>\nu +mL+n$, treating the remainder $q'_J$ as last
described, and the homogeneous terms with $j<J$ as
described above. 

This shows how (1.16) is obtained in general, reduced to
(1.19) in case (4) (and (3)). Clearly, each $\tilde c_j(x)$ depends only
on the strictly homogeneous term of degree $-mN-j$ in $q$,  hence on the
strictly homogeneous terms of the first $j+1$ 
orders in $A$ and $P$.

It remains to show the formulas 
(1.17), (1.20) for $\tilde c''_0(x)$. Here we go back to the
splitting (1.25), applied to 
the complete symbol $q(x,\xi ,\mu )$, that we examine with special care.
It is observed in \cite{GS1, p.\ 501} that when $P$ has symbol $p(x,\xi )$, the
symbol $\tilde q(x,\xi ,\mu )$ of $(P+\mu ^{m})^{-1}$ is a sum $(p_m(x,\xi
)+\mu ^m )^{-1}+\tilde q'$, where $(p_m+\mu ^m )^{-1}\in S^{-m,0}\cap S^{0,-m}$ and $\tilde q'$ has not only lower order but also
better decrease in $\mu $: $\tilde q'\in S^{-m-1,0}\cap
S^{m-1,-2m}$ (since it is constructed from terms containing
at least two powers of $(p_m+\mu ^m )^{-1}$). A similar phenomenon
holds for the $N$'th power of $(P+\mu ^m)^{-1}$; its symbol is 
a sum $(p_m+\mu ^m )^{-N}+\tilde q_{(N)}'$, where $(p_m+\mu
^m )^{-N}\in S^{-mN,0}\cap S^{0,-mN}$ and $\tilde q'_{(N)}\in 
S^{-mN-1,0}\cap 
S^{m-1,-m(N+1)}$. The composition of $A$ with
$\operatorname{OP}(\tilde q'_{(N)})$ gives an operator with symbol in
$S^{\nu -mN-1,0}\cap 
S^{\nu +m-1,-m(N+1)}$ (depending on $\mu ^m$ rather than $\mu $); by
the preceding analysis, it has a diagonal kernel
expansion of the form
$$
K(x,x ,\mu )
\sim
\sum_{ j\in \Bbb N  }  d_{ j}(x)\mu  ^{\nu -1+n -j-mN}+ 
\sum_{k\ge 1}\bigl( d'_{ k}(x)\log \mu  +d''_{
k}(x)\bigr)\mu ^{-mk-mN},
\tag1.30
$$
where the sum over $k$ starts with $k=1$, so that
$\mu ^{-mN}$ appears at most with a local coefficient from the series
over $j$. In the
consideration of the symbol composition 
$$
\multline
a(x,\xi )\circ(p_m(x,\xi)+\mu ^m )^{-N}=
a(x,\xi )(p_m(x,\xi)+\mu ^m )^{-N}\\
+
\sum _{1\le |\alpha |<M}\tfrac {(-i)^{|\alpha |}}{\alpha !}\partial
_\xi ^\alpha a(x,\xi )\partial _x^\alpha (p_m(x,\xi)+\mu ^m )^{-N}
+r_M(x,\xi ,\mu ),
\endmultline\tag1.31
$$ we observe that whenever a 
derivative (in $x$ or $\xi $) hits $(p_m+\mu ^m )^{-N}$, its $d$-index
is lowered (since the resulting expression contains at least one more
power of    $(p_m+\mu ^m )^{-1}$) --- and the same is true for
the remainder $r_M$
(constructed by Taylor expansion as in standard proofs of the composition
rule). So again, the part  $\sum _{1\le |\alpha |<M}\tfrac {(-i)^{|\alpha |}}{\alpha !}\partial
_\xi ^\alpha a\partial _x^\alpha (p_m+\mu ^m )^{-N}
+r_M
$ gives a kernel expansion of the form (1.30),
where $\mu ^{-mN}$ appears at most in the series over $j$, with a
local coefficient. (One could avoid this step by taking the symbol of
$(P+\mu ^m)^{-N}$ in $y$-form, found from the conjugate transpose of
the symbol of $(P^*+\mu ^m)^{-N}$.)

It remains to consider 
$\operatorname{OP}(a(x,\xi )(p_m(x,\xi )+\mu ^m
)^{-N})$. Here
we remark that$$
\aligned
\tilde q''&\equiv\mu ^{-mN}-(p_m+\mu ^m )^{-N}=\big((p_m+\mu ^m)^N-\mu ^{mN})(p_m+\mu ^m
)^{-N}\mu ^{-mN}\\
&=\sum_{1\le l\le N}\tbinom N l p_m^{l}\mu ^{m(N-l)}(p_m+\mu ^m
)^{-N}\mu ^{-mN};
\endaligned \tag1.32
$$ 
a sum of terms in
$S^{m(l-N), -ml}\cap S^{ml,-m(l+N)}$, $l=1,\dots,N$; the sum is in 
$S^{0, -m}\cap S^{mN,-m(1+N)}$. 
Now write$$
\gathered a(x,\xi )=a'(x,\xi )+a''(x,\xi ),\\
a'=\sum_{j\le \nu +n}a_{\nu -j},\quad a''=a-\sum_{j\le \nu
+n}a_{\nu -j}.
\endgathered \tag 1.33
$$
Here $a''$ is of order $\nu -J<-n$ ($J=\max\{n+[\nu ]+1,0\}$) and
defines a trace-class operator. For 
$\operatorname{OP}(a''(p_m+\mu ^m)^{-N})$, the diagonal kernel is $$
\int_{\Bbb R^n}a''(p_m+\mu ^m)^{-N}\,\d\xi =\mu ^{-mN}\int_{\Bbb
R^n}a''\,\d\xi -\int_{\Bbb R^n}a''\tilde q''\,\d\xi. 
$$
Since $a'\tilde q''\in S^{\nu -J, -m}\cap S^{\nu -J+mN,-
m(1+N)}$, the last integral
gives a series as in (1.30), now with $\nu -1$ replaced by $\nu -J<-n$
so that there is no term with
$\mu ^{-mN}$. Thus the contribution 
from $a''$ to the coefficient
of $\mu ^{-mN}$ is $\int a''\,\d\xi $, the last
parenthesis in (1.18).

For $a'(p_m+\mu ^m)^{-N}$, we know that the integral over $|\xi |\ge
\mu $ of each \linebreak$a_{\nu -j}(p_m+\mu ^m)^{-N}$ gives a local
term, as in 
the sum over $j$ in (1.16). For the integral over
$|\xi |\le \mu $, we consider the two parts $\mu ^{-mN}a'$ and
$-a'\tilde q''$ separately. The latter gives a sum of expressions
$$
\tbinom Nl \mu ^{-ml}\int_{|\xi |\le \mu }a'p_m^l(p_m+\mu
^m)^{-N}\,\d\xi ,\tag1.34
$$
where we find as in the analysis of $q$ above that the integral alone
produces terms as in (1.26) and (1.28), where nonlocal contributions
start at the power $\mu ^{-mN}$. Thanks to the extra factor $\mu
^{-ml}$ in front ($l\ge 1$), (1.34) on the whole contributes locally
(as in the sum over $j$ in (1.30)) to the coefficient of $\mu ^{-mN}$.

Finally, we study $\mu ^{-mN}\int_{|\xi |\le \mu }a'(x,\xi
)\,\d\xi $. The integral over $|\xi |\le 1$ simply gives 
 $\mu ^{-mN}\int_{|\xi |\le 1 }a'(x,\xi
)\,\d\xi $.
The integral over $1\le |\xi |\le \mu $ of each
homogeneous term is analyzed as in
(1.28); here the contribution from $a_{\nu -j}$ to the coefficient of
$\mu ^{-mN}$ is $-\frac 1{\nu -j+n}\int_{|\xi |=1}a_{\nu -j}(x,\xi
)\,\d S(\xi )$ if $\nu -j\ne -n$, zero if $\nu -j=-n$.

Adding the contributions from $a'$ and $a''$ we find (1.17) with
(1.18). This completes the proof of (i) and (ii).

For (iii), we just have to check the local contributions found along
the way in the preceding considerations. In case (1), the series in
$j$ (as in (1.30)) begin below the power $\mu ^{-mN}$, and in case (2), they
contain only noninteger powers of $\mu $. In the cases (3)
and (4), it is checked as in the beginning of the proof that the
series in $j$ contain only odd powers of $\mu $. So in all these
cases, the only contributions to the coefficient of $\mu ^{-mN}$ come
from $\tslint a(x,\xi )\,\d\xi $. 
\qed
\enddemo

Point (iii) in cases (1) and (2) was shown already by Lesch in \cite{L};
he introduced the notation
$\tslint a(x,\xi )\,\d\xi $ with the following brief description: It
equals the $\mu $-independent
term $\tilde a''_0(x)$ in the asymptotic expansion for $\mu
\to\infty $:$$ 
\int_{|\xi |\le \mu }a(x,\xi )\,\d\xi \sim \sum_{j\in\Bbb N, j\ne \nu
+n}\tilde a_j(x)\mu ^{\nu +n-j}+\tilde a_0'(x)\log\mu +\tilde a''_0(x).\tag1.35
$$
This description is clearly consistent with the
calculation of $\tslint a\,\d\xi $ in the above proof. Also the
notation $\operatorname{LIM}_{\mu \to\infty }\int_{|\xi |\le \mu
}a(x,\xi )\,\d\xi $ is used. The concept is related to  Hadamard's
definition of the finite part --- {\it partie finie} --- of certain
integrals \cite{H, p\. 184 ff.}. 

Lesch moreover shows that in the cases (1)
and (2), the density 
$\tslint a(x,\xi 
)\,\d\xi |dx|$ associated with $A$ is invariant under coordinate changes. 
In fact, he shows this also when $A$ is given by an amplitude
function $a(x,y,\xi
)$ (a symbol in
$(x,y)$-form); then $\tslint a(x,x,\xi )\,\d\xi |dx|$ is invariant.
The proof of the invariance extends to the cases
(3) and (4), since in the proof of \cite{L, Lemma 5.3}, the logarithmic
contributions (the sum over $l$) in Prop\. 5.2 vanish because the
terms of order $-n$ are odd in $\xi $.
So in all the cases,
$A $ defines a density
$\omega_{\operatorname{TR}}(A)$
described in local coordinates by $$
\omega_{\operatorname{TR}}(A)=\slint a(x,\xi )\,\d\xi |dx|\text{
resp\. }\omega _{\operatorname{TR}}(A)=\slint a(x,x,\xi )\,\d\xi\, |dx|,\tag 1.36 
$$
when $A$ has symbol in $x$-form $a(x,\xi )$, resp\. in  $(x,y)$-form
$a(x,y,\xi )$.

Note that in (1.18) and (1.33), one can replace the sum over $j\le \nu +n$ by
the sum over $j\le J$ for any choice of $J\ge \nu +n$, since$$
\slint a_{\nu -j}\,\d\xi =\int a_{\nu -j}\,\d\xi =\int_{|\xi |\le
1}a_{\nu -j}\,\d\xi +\tfrac 1{j-\nu -n}\int_{|\xi |=1}a_{\nu -j}\,\d
S(\xi )
$$
for $j>\nu +n$, by integrability and homogeneity.
Note also that $\tslint  a(x,\xi )\,\d\xi =0 $ when $a$ is polynomial
in $\xi $ (reconfirming the fact from \cite{GS1, Th.\ 2.7} that the
coefficient of $(-\lambda )^{-N}$ is local when $A$ is a differential
operator). 

For convenience, we recall that when $A=I$, the
coefficients of $(-\lambda )^{\frac {n-j}m-N}$ for $j<m+n$ are
simply$$
\tilde c_j(x)=\int_{\Bbb R^n}q^h_{-mN-j}(x,\xi ,1)\,\d\xi , \quad
0\le j<m+n,\tag 1.37
$$
where $q^h_{-mN-j}$ is the strictly homogeneous version of
$q_{-mN-j}$; a direct proof goes as in e.g.\ \cite{G1, Th.\
3.3.5, cf\. (3.3.31), (3.3.39)} or \cite{GS1, (2.16)}, using that the
$q^h_{-mN-j}$ are 
integrable at $\xi =0$ for $j<m+n$. (Also noninteger $m>0$ are allowed
here.) This includes $\tilde c_n(x)$, the coefficient of $(-\lambda
)^{-N}$; here $C_0(I,P)=\int \tr \tilde c_n(x)\,dx$.

\example{Remark 1.4} The formulas (1.17), (1.18), (1.20) in
Theorem 1.3 extend to the 
case where $A$ is given by a symbol $a(x,y,\xi )$ in   $(x,y)$-form;
then  $a(x,\xi )$ in the formulas is replaced by $a(x,x,\xi )$. To
see this, we just 
need a supplement to the last part of the above proof. 
We split in two parts $$
a(x,y,\xi )=a'(x,y,\xi )+a''(x,y,\xi ),\quad a'=\sum_{j\le \nu
+n}a_{\nu -j},\tag1.38
$$ as in (1.33). The trace-class part is easily dealt with:
As in (1.32), $\mu ^{-mN}-(P+\mu ^{m})^{-N}$ has symbol in
$S^{0,-(m+1)N}$, so its composition with
$\operatorname{OP}(a''(x,y,\xi ))$ has symbol in $S^{\nu
-J,-(m+1)N}$ with $\nu -J<-n$ (by the general composition rules for these
symbol spaces), hence gives a diagonal kernel expansion as in (1.30)
with $\nu -1$ replaced by $\nu -J$; it contains no term with $\mu
^{-mN}$. The diagonal kernel of $\operatorname{OP}(a''(x,y,\xi ))\mu
^{-mN}$ is simply $\mu ^{-mN}\int a''(x,x,\xi )\,\d\xi $.

For $\operatorname{OP}(a')$, we
observe that by the rules of calculus for $\psi $do's, $
\operatorname{OP}(a'(x,y,\xi ))=\operatorname{OP}(a_1(x,\xi ))+\operatorname{OP}(a_2(x,y,\xi ))$,
where $$
\aligned
a_1(x,\xi )&=\sum_{|\alpha |<M}\tfrac{(-i)^{|\alpha |}}{\alpha !}\partial
_\xi ^\alpha \partial _x^\alpha a'(x,y,\xi )|_{y=x},\\
a_2(x,y,\xi )&=\sum_{|\alpha |=M}\tfrac{M(-i)^{|\alpha |}}{\alpha
!}\partial _\xi ^\alpha \int_0^1(1-t)^{M-1}\partial_y^\alpha
a'(x,x+t(y-x),\xi )\,dt; 
\endaligned\tag 1.39
$$
we take $M> \nu +n+1$.
We can take
the symbol $\tilde q$ of $(P+\mu ^m)^{-N}$ in $y$-form.
The first term in the first line of (1.39) is $a'(x,x,\xi )$, a symbol in
$x$-form, whose effect is
as described in the theorem; this gives the value in (1.17)
with $a(x,\xi )$ replaced by $a'(x,x,\xi )$.
The other terms in the first line 
are also in $x$-form, now with a power $\partial _\xi ^\alpha $ in front. When
the corresponding operators are composed with
$\operatorname{OP}(\tilde q(y,\xi ,\mu ))$ and the kernel is
calculated, we can perform an 
integration by parts w.r.t\. $\xi $, placing the derivative on
$\tilde q$.
As noted earlier, the derivatives of $\tilde q$ have symbols with $d$-index
$\le -(m+1)N$, so the resulting integrals have expansions as in (1.30),
giving only local contributions to the
coefficient of $\mu ^{-mN}$. They vanish in the cases (1)--(4).

As for $a_2$, it is in $(x,y)$-form and equals a sum of $\xi
$-derivatives  \linebreak$\sum_{i=1}^n\partial _{\xi
_i}a_{2,i}(x,y,\xi )$ with $a_{2,i}$ of order $<-n$. Here the
considerations on $a''(x,y,\xi )$ apply, and moreover, the
contributions to (1.18) vanish since $\int_{\Bbb
R^n}\partial _{\xi _i}a_{2,i}(x,x,\xi )\,\d\xi =0$ as the integral of
a derivative.

\endexample

\proclaim{Corollary 1.5} Consider operators $A$ and $P$ on the
manifold $X$, then there is an asymptotic expansion {\rm (1.1)} for
$N>(n+\nu )/m$. 
In the localized situation,
$$
C_0(A,P)=\tilde c_{\nu +n}+\tilde c''_0=\int_{\Bbb
R^n}\slint \tr a(x,\xi )\,\d\xi dx +\text{ local terms};\tag1.40
$$
the local terms depend only on the first $n+[\nu] +1$ strictly homogeneous
terms in the symbols of $A$ and $P$. Here $\tr$ denotes fiber trace.
(If the symbol of $A$ is
in $(x,y)$-form, the formula holds with $a(x,x,\xi )$ instead.)

{\rm (ii)} If, moreover,  $A$ and $P$ are as in Definition
{\rm 1.1} 
with {\rm (3)} or {\rm (4)}, then the expansion {\rm (1.1)} reduces
to 
the form$$
\Tr(A(P-\lambda )^{-N})
\sim
\sum_{ j\in \Bbb N ,\,j-n-\nu \text{ odd}} \tilde c_{ j}(-\lambda )
^{\frac{\nu +n -j}m-N}+  
\sum_{k\in \Bbb N}\tilde c''_{
k}(-\lambda ) ^{{ -k}-N}.
\tag1.41
$$
In particular, $\zeta (A,P,s)$ has no poles at integers $s$.

{\rm (iii)} In each of the cases {\rm (1)--(4)} in Definition {\rm
1.1}, $\tilde c'_0=\tilde c_{\nu +n}=0$ and 
(cf\. {\rm (1.36)})$$
C_0(A,P)=\tilde c_0''=\int _X\tr \omega _{\operatorname{TR}}(A).\tag 1.42
$$

\endproclaim 

\demo{Proof} Since $A(P-\lambda )^{-N}$ is trace-class,
$\Tr(A(P-\lambda )^{-N})$ can be expressed in the chosen local
coordinates as the 
integral in $x$ of the fiber trace of the kernel diagonal value. Then
the corollary follows 
directly from Theorem 1.3 by integration in $x$. For the formula
(1.42) we use the information leading to (1.36).\qed
\enddemo 

Note in particular that we have obtained that $\tilde c''_0$ depends
only on $A$ (not on the auxiliary operator $P$) in the cases
(1)--(4). See also Remark 1.8 further below.

\example{Remark 1.6} Note that in all the cases
(1)--(4) in Definition 1.1, $\zeta (A,P,s)$ is regular at $s=0$, and$$
\zeta (A,P,0)= c''_0-\Tr(A\Pi _0(P))=\operatorname{TR}A-\Tr(A\Pi
_0(P)).\tag1.43 
$$
Case (3) plays an important role in \cite{O3}. The operators $A$ in case
(4) do not in themselves form an algebra (neither do the operators in
the cases (1) and (2)), but we think that the
definition is of interest anyway. As examples of
case (4) we mention eta functions on an even-dimensional
manifold $X$: Let $D$ be a selfadjoint first-order elliptic
differential operator on $X$, then the eta function $\eta (D,s)$ is
defined for $\operatorname{Re}s$ large by$$
\eta (D,s)=\Tr(D|D|^{-s-1})=\Tr(D(D^2)^{-\frac12}(D^2)^{-\frac
s2})=\zeta (D|D|^{-1},D^2,\tfrac s2);
$$ 
here $D|D|^{-1}$ (defined to vanish on the nullspace of $D$) is of
order 0 and even-odd.
For the meromorphic extension according to (1.2), the locally
determined coefficients $c_{n}$ and $c'_0$ vanish due to parity, and
$D\Pi _0(D^2)=0$, so $$
\eta (D,0)=c''_0=\operatorname{TR}(D|D|^{-1}).\tag1.44
$$
\endexample

It will now be shown that the expression $\operatorname{TR}A$ defined in
Definition 1.1 vanishes on commutators:

\proclaim{Theorem 1.7} Let $A $ and $A'$ be classical $\psi $do's
of orders $\nu  $ resp\. $\nu '\in\Bbb R$, and let $P$ be a classical
elliptic $\psi $do 
of even order $m>0$ such that the principal symbol has no eigenvalues
on $\Bbb R_-$.
Then $$
\operatorname{TR}([A ,A'])=0\tag1.45
$$
holds in the following cases:
\roster
\item "(1$'$)" $\nu  +\nu '<-n$ ($\nu  +\nu '<1-n$ if the principal symbols commute). 
\item "(2$'$)" $\nu  +\nu '\notin \Bbb Z$.
\item "(3$'$)" $\nu$ and $\nu '\in\Bbb Z$,
$A $ and $A'$ are 
both even-even or both even-odd, and $n$ is odd.
\item "(4$'$)" $\nu $ and $\nu '\in\Bbb Z$,
 $A $ is 
even-odd, $A'$ is even-even, and $n$ is even.
\endroster
\endproclaim 

\demo{Proof} The case (1$'$) is an immediate consequence of the definition,
since   $\operatorname{TR}([A ,A'])=\Tr ([A ,A'])$ and the standard
trace vanishes on commutators. For the other
cases, we rewrite by use of suitable resolvent formulas and
cyclic permutation, taking  $P$ even-even in case (3) and (4):
In view of the identity $$
\partial _\lambda ^k(P-\lambda )^{-1}=k!(P-\lambda )^{-k-1},\tag 1.46
$$we have that$$\aligned
\Tr\bigl([A ,A'](P-\lambda )^{-N}\bigr)
&=\tfrac1{(N-1)!}\Tr\bigl(
\partial _\lambda ^{N-1}A A'(P-\lambda )^{-1}-\partial _\lambda
^{N-1}A 
(P-\lambda )^{-1}A'\bigr)\\
&=\tfrac1{(N-1)!}\Tr \bigl(
A \partial _\lambda ^{N-1}((P-\lambda )^{-1}[P,A'](P-\lambda
)^{-1})\bigr)\\
&=\Tr\bigl(A \sum_{0\le M<N}c_{MN}(P-\lambda
)^{-1-M}[P,A'](P-\lambda )^{-N+M} \bigr).
\endaligned\tag1.47$$
Let$$
Q(\mu )=A \sum_{0\le M<N}c_{MN}(P+\mu ^m)^{-1-M}[P,A'](P+\mu
^m)^{-N+M}.\tag1.48 
$$
By the rules of calculus in \cite{GS1}, 
$Q(\mu )$ has symbol in $$
S^{\nu  +\nu '-mN,0}\cap S^{\nu  +\nu '+m,-m(N+1)}.\tag1.49
$$
Here $\nu +\nu '$ can be replaced by $\nu +\nu '-1$ if $P$ is scalar,
but this is without importance for the special trace term we are
investigating. 
Since $Q$ actually only depends on $\lambda =-\mu ^m$, the symbol
expansion like that of $f'$ in (1.15) has only powers that are
multiples of $m$. 
The important thing here is that the lowest $d$-index, $d=-m(N+1)$, is
lower than
that of $A(P+\mu ^m)^{-N}$ itself. 
By \cite{GS1, Th\. 2.1}, there is 
then a trace expansion
$$\aligned
\Tr\bigl( [A ,A']&(P-\lambda )^{-N}\bigr)=\Tr Q(\mu )\\
&\sim
\sum_{ j\in \Bbb N  } \tilde b_{ j}(-\lambda ) ^{\frac{\nu  +\nu '+n -j}m-N}+ 
\sum_{k\ge 1}\bigl( \tilde b'_{ k}\log (-\lambda ) +\tilde b''_{
k}\bigr)(-\lambda ) ^{{ -k}-N}.
\endaligned\tag1.50
$$

When $\nu  +\nu '\notin \Bbb Z$, there is no term of the form
$c(-\lambda )^{-N}$, so $\operatorname{TR}([A ,A'])$ vanishes
according to our definition. This takes care of the case (2$'$).

In the cases (3$'$) and (4$'$), we note that
when $A $ 
and $A'$ have the same alternating parity, then $A A'$ and $[A ,A']$
are even-even, whereas
when they have opposite parities, $A A'$ and $[A ,A']$
are even-odd.
Then we can use the information from Theorem 1.3 that the terms in
the series over $j$ vanish for $\nu +\nu '+n-j$ even; this holds in
particular for the constant term, where $j=\nu +\nu '+n$.
\qed
\enddemo 

The result of Theorem 1.7 is essentially known from \cite{KV} in the
cases (1$'$)--(3$'$), but our proof is different; it will be
generalized to log-polyhomogeneous operators in Section 3. 

We note in passing that the proof that $\operatorname{TR}A$ in
Definition 1.1 is independent of the choice of
$P$ could also be based on resolvent rules instead of the
painstaking analysis of its value:
Since $(P-\lambda )^{-1}-(P'-\lambda )^{-1}=(P-\lambda )^{-1}(P'-P)(P'-\lambda )^{-1}$,
we can write
$$\aligned
\Tr\bigl(A&(P-\lambda )^{-N}\bigr)-\Tr\bigl(A(P'-\lambda )^{-N}\bigr)\\
&=\tfrac1{(N-1)!}\Tr\bigl(
A\partial _\lambda ^{N-1}((P-\lambda )^{-1}(P'-P)(P'-\lambda )^{-1})\bigr).
\endaligned\tag1.51
$$
This operator family has symbol in 
$S^{\nu -mN,0}\cap S^{\nu +m,-m(N+1)}$,
hence has a trace expansion as in (1.50) with $\nu +\nu '$ replaced
by $\nu $.
From this, one can reason exactly as in the proof of Theorem 1.7.

\example{Remark 1.8} In the above considerations, we have kept the
order of $P$ fixed, 
equal to an even number. One can in fact show trace expansions with a
similar structure as in (1.1), (1.2) when the order of $P$ is an
arbitrary $m\in \Bbb R_+$, 
cf\. Loya \cite{Lo}, 
Grubb and Hansen \cite{GH}, which could have been taken as the point of
departure. On the other hand, when $P$ is of an arbitrary order $m >0$,
$\zeta (A,P,s)=\zeta (A, P^{ {m'}/m }, s')$ for $s'=sm/{m'}$, so by a scaling of the complex variable $s$ one can reduce to a
situation  where the order of $P$ is a given even number, at least
for positive selfadjoint operators.
\endexample

\bigskip
\head 2. A quasi-trace
\endhead
\medskip
The functional $\operatorname{TR}$ does not extend to the general case $\nu
\in\Bbb Z$ as a 
trace (cf\. \cite{KV}, Lesch \cite{L}). Yet it is still possible to make
some further observations on the integer
order case. Consider $C_0(A,P)$ defined in (1.4).
It coincides
with $\operatorname{TR}A$ in the cases in Definition 1.1, 
 but depends in general on $P$.
However, $C_0(A,P)$ has the independence of $P$ and the commutator
property in a
weaker sense, namely:  

\proclaim{Proposition 2.1}

{\rm (i)} Let $A$ be a classical $\psi $do
of order $\nu \in\Bbb
Z$, and let $P$ and $P'$ be classical elliptic $\psi $do
of positive orders $m$ and $m'$ such that the principal symbols have no eigenvalues
on $\Bbb R_-$. Then $C_0(A,P)-C_0(A,P')$ is 
locally determined. More precisely, it depends solely on the strictly homogeneous parts
of the first $\nu +n+1$ homogeneous
terms in each of the 
symbols of $A$, $P$ and $P'$; it vanishes if $\nu <-n$.

{\rm (ii)} Let $A $ and $A'$ be classical $\psi $do's
of order $\nu  $ resp\. $\nu '$ such that  $\nu  +\nu '\in\Bbb
Z$, and let $P$ be a classical elliptic $\psi $do
of even order $m>0$ such that the principal symbol has no eigenvalues
on $\Bbb R_-$.
Then $C_0([A ,A'],P)$ is
locally determined. More precisely it 
depends solely on the strictly homogeneous parts of the first $n+\nu
+\nu '+1$ homogeneous terms in each of the
symbols of $A $, $A'$ and $P$; it vanishes if $\nu +\nu '<-n$.

\endproclaim 

\demo{Proof} 
(i). By Remark 1.8, we may assume that $P$ and $P'$ have the same
even order $m$. Then the statement follows directly  from Corollary 1.5 (i). 

(ii). As noted in the proof of Theorem {\rm 1.7},
$\Tr([A ,A'](P-\lambda )^{-N})$ has an expansion (1.50). Since $\nu
+\nu '$ is 
integer, and the sum over $k$ begins with $k=1$,$$
C_0([A ,A'],P)=\tilde b_{\nu  +\nu '+n},\tag2.1
$$
which is locally determined as stated.
\qed\enddemo 

Note that the expressions $C_0(A,P)-C_0(A,P')$ and $C_0([A,A'])$ are
{\it pointwise} locally determined, in the sense that they can be
calculated as integrals in $x$ of locally determined functions (in
local coordinates). 

The proposition shows that $C_0(A,P)$ is in general somewhat like a
trace, just {\it modulo local contributions}.  
The values of $C_0([A ,A'],P)$ and
$C_0(A,P)-C_0(A,P')$ can be described in terms of
certain residues, as 
we shall recall in Proposition 3.1 below.
However, for the sake of more general situations, we believe that it
has an interest to introduce a notion of
quasi-trace as follows:

\proclaim{Definition 2.2} Let $X$ be an $n$-dimensional compact
manifold without boundary, provided with a $C^\infty $ vector bundle
$E$, let $A$ run through an algebra of $\psi
$do's of orders $\nu $ on the sections of $E$, and let $P$ run through an auxiliary
family of 
elliptic $\psi $do's in $E$ without
principal symbol eigenvalues on $\Bbb R_-$.
Consider a function 
$f(A,P)$ such that for each $P$, it is a linear functional on the
$A$'s. We say 
that $f$ is a quasi-trace if {\rm (i)} and {\rm (ii)} hold:
\roster
\item "{\rm (i)}" $f(A,P)-f(A,P')$ depends only on the strictly
homogeneous symbols of the first $[\nu ]+n+1$ degrees in the
symbols of $A$, $P$ and $P'$, and vanishes for $\nu <-n$.
\item "{\rm (ii)}" $f([A ,A'],P)$ depends only on the strictly
homogeneous symbols of the first  $[\nu +\nu ']+n+1$ degrees in the
symbols of $A$, $A'$ and $P$, and vanishes for $\nu +\nu '<-n$.
\endroster

\endproclaim

By the preceding results, $C_0(A,P)$ is a quasi-trace in this sense.
(In the
formulation of (i) and (ii), when $j\le
0$, the set of symbols of the first $j$ degrees is understood to be
empty. Thus since $f$ is linear in $A$, the statement on the
vanishing is a consequence of the preceding statement.)

As shown in (i) of Theorem 1.3, $C_0(A,P)$ moreover has a pointwise
description, where it can be obtained, modulo
local contributions, as an integral in $x$
of the function   $\tslint \tr a(x,\xi )\,\d\xi $, defined from the symbol
$a(x,\xi )$ of $A$ in the chosen local
coordinates.

Scott in \cite{Sco}
uses the name $\operatorname{TR}_\Delta $ for a concept like
$C_0(A,P)$ with
$P=\Delta $. We have recently been informed that the constant $C_0(A,P)$
(and its generalizations to $b$-calculi)
plays an important role in the manuscript of
Melrose and Nistor \cite{MN} and subsequent works, where it is called a
{\it regularized trace}, denoted $\widehat{\operatorname{Tr}}(A)$ or
$\Tr_P(A)$.  It is taken up in a physics context under the name of a
weighted trace $\tr ^P(A)$ in Cardona, Ducourtioux, Magnot, Paycha
\cite{CDMP}, \cite{CDP}.

The concept can also be useful
when one has a vanishing property of the relevant local
contributions, as  e.g\. in \cite{G3}.
A generalization to manifolds with boundary is worked out in Grubb
and Schrohe \cite{GSc}.
\bigskip
\head 3. Higher Laurent coefficients of zeta functions,
log-po\-ly\-ho\-mo\-ge\-ne\-ous symbols
\endhead
\medskip
Throughout the following, we make the extra assumption that $P$ {\it is
invertible}. (This makes the 
statements simpler. In general, one can replace $P$ by $P+P_0$, where
$P_0$ is the operator with kernel $K(P_0,x,y)=\sum_{1\le j\le \dim\ker
P}\varphi '_j(x)\varphi ^*_j(y)$, where the $\varphi _j'$ and
$\varphi _j$ denote the zero eigensections of $P^*$ resp.\ $P$, and
correct for terms 
stemming from $P_0$ afterwards.)
Then $\zeta (A,P,s)=\Tr(AP^{-s})$ has the
Laurent expansion (1.5) at $s=0$.
The noncommutative residue $\operatorname{res}(A)=m\cdot C_{-1}(A,P)$
and the quasi-trace 
$C_0(A,P)$ have been discussed above, but 
also the next coefficient $C_1(A,P)$ is of particular interest; in
the case $A=I$
it equals minus the so-called zeta-determinant of $P$:$$
\log\det P=- \partial _s\zeta (I,P,0)=-C_1(I,P).\tag 3.1
$$
Not only 
this coefficient, but the whole Laurent expansion (1.5) can be
described by use of 
functional calculus combined with the work of Lesch \cite{L}, as we
show in the following.

According to Seeley \cite{S}, the complex powers of $P$ are defined 
by$$
\aligned
P^{-s}&=\tfrac{ i}{2\pi
}\int_{{\Cal C}}\lambda 
^{-s}(P-\lambda )^{-1} \,d\lambda,\text{ for $\operatorname{Re}s>0$ 
},\\
P^{-s}&=P^{-s-N}P^N=P^NP^{-s-N}\text{ in general};
\endaligned
\tag3.2
$$ 
here $\Cal C$ is a curve in $\Bbb C\setminus \overline{\Bbb R}_- $
encircling the spectrum of $P$ in the positive direction (replace
intervals of $\crm $ by small half-circles around the finitely many possible
eigenvalues of $P$ on $\Bbb R_-$).
One defines $\log P=-\partial _sP^{-s}|_{s=0}$ (cf\. e.g\. \cite{KV},
\cite{O1}); then $$ 
\aligned
\partial _sP^{-s}&=-\log P\, P^{-s},\text{ where}\\
\log P\, P^{-s}&=\tfrac{ i}{2\pi
}\int_{{\Cal C}}\log \lambda \,\lambda 
^{-s}(P-\lambda )^{-1} \,d\lambda\text{ when $\operatorname{Re}s>0$ 
}.
\endaligned
\tag3.3
$$
It is known (cf\. e.g\. \cite{O1}) that $\log P$ is
a $\psi $do such that in local coordinates,$$
\text{symbol of }\log P\sim m\log[\xi ]I +\sum_{j\ge 0}b_{-j}(x,\xi ) 
\tag3.4$$ 
with $b_{-j}$ homogeneous in $\xi $ of degree $-j$, 
$[\xi ]$ smooth positive and equal to $|\xi |$ for $|\xi |\ge 1$. These
homogeneous terms are derived straightforwardly from the homogeneous
terms in the symbol of $(P-\lambda )^{-1}$; in particular, $b_{-j}$ 
is determined from the first $j+1$ homogeneous terms in the symbol of $P$. 

We remark, for the convenience of the reader, that one can show the
following precisions of the 
quasi-trace property of $C_0(A,P)$, using $\log P$ and the operator
family $P^{-s}$:

\proclaim{Proposition 3.1} Let $A$, $A'$, $P$ and $P'$ be as in
Proposition {\rm 2.1}. Then
$$\gather
C_0(A,P)-C_0(A,P')=\tfrac 1m \operatorname{res}(A(-\log P+\log P')),\tag3.5
\\
C_0([A,A'],P)=\tfrac 1m \operatorname{res}(A[\log P,A']).\tag3.6
\endgather$$
\endproclaim

\demo{Proof} (3.5) is shown in \cite{O1} and \cite{KV}. (3.6)
is shown in \cite{MN}, where its generalizations to the $b$-calculus play
an important role; related formulas have appeared in the physics litterature
(see e.g\. Mickelsson \cite{M}, Cederwall, Ferretti, Nilsson and
Westerberg \cite{CFNW}). 
The following method of proof is deduced from \cite{MN}. 

We first observe the following consequence of \cite{Gu} and \cite{W} (in
particular \cite{Gu, Th.\ 7.1}): Let $B(s)$ be a holomorphic family of
classical $\psi $do's of order $\alpha  -s$ (for some real $\alpha
$), with $B(0)=0$. Then the
meromorphic extension of $\Tr B(s)$ (holomorphic for
$\operatorname{Re}s> n+\alpha $, the extension again denoted $\Tr
B(s)$), satisfies:$$
\lim_{s\to 0}\Tr B(s)=\operatorname{res}B'(0).\tag 3.7
$$
For, setting $C(s)=\frac1s (B(s)-B(0))=\frac1s B(s)$, we have that
$C(s)$ is a holomorphic family of classical $\psi $do's of order
$\alpha -s$ with $C(0)=B'(0)$. Moreover,  $Q(s)=\frac1s
(C(s)-C(0)P_1^{-s})$ is 
another holomorphic family of classical $\psi $do's of order $\alpha -s$,
when $P_1$ is taken as an elliptic positive selfadjoint $\psi $do of
order 1. Then$$
B(s)=sC(s)=sC(0)P_1^{-s}+s(C(s)-C(0)P_1^{-s})=sB'(0)P_1^{-s}+s^2Q(s).
$$
Now$$
\lim_{s\to 0}\Tr (sB'(0)P_1^{-s})=\operatorname{res}B'(0),
$$
by definition of the noncommutative residue, and $$
\lim_{s\to 0}\Tr (s^2Q(s))=0,\tag3.8
$$
by \cite{Gu, Th.\ 7.1} (assuring that $\Tr Q(s)$ extends
meromorphically to $\Bbb C$, with {\it simple} poles lying in $\alpha -\Bbb
Z$). This shows (3.7).

Applying (3.7) to the holomorphic family $A(P^{-s/m}-(P')^{-s/m})$ of
order $\nu -s$, we see that$$
\aligned
\lim_{s\to 0}\Tr (A(P^{-s}-(P')^{-s}))&=
\lim_{s\to 0}\Tr (A(P^{-s/m}-(P')^{-s/m}))\\
&=\tfrac
1m\operatorname{res}(A(-\log P+\log P')),\endaligned
$$
since $\partial _sP^{-s/m}=-\frac1m\log P\, P^{-s/m}$ (note that $-\log P+\log
P'$ is classical in view of (3.4)); this shows (3.5).

For (3.6), we apply (3.7) to the family $A[A',P^{-s/m}]$ of order
$\nu +\nu '-s$. Noting that for $\operatorname{Re}s$ large, $$
\Tr([A,A']P^{-s})=\Tr(AA'P^{-s})-\Tr(AP^{-s}A')=\Tr(A[A',P^{-s}])
$$
by cyclic permutation, we find:$$
\lim_{s\to 0}\Tr ([A,A']P^{-s})=\lim_{s\to 0}\Tr (A[A',P^{-s/m}])=
\tfrac1m \operatorname{res}(A[\log P,A']),
$$
which shows (3.6). Here $[\log P, A']$ is classical in view of (3.4).\qed
\enddemo 

The meromorphic extension $\Tr A(s)$ for a holomorphic family $A(s)$
of order $\alpha -s$ coincides with $\TR A(s)$ when $\alpha
-s\notin\Bbb Z$, by \cite{KV, Th.\ 3.1}.

[Added November 2005: A proof of Proposition 3.1 based directly on
resolvent information is given in \cite{G4}.]

Consider now the higher derivatives of $P^{-s}$:
$$\aligned
\partial _s^2 P^{-s}&=(-\log P)^2 P^{-s},\\
&\vdots\\
\partial _s^{l} P^{-s}&=(-\log P)^{l}P^{-s},
\endaligned \tag3.9
$$
which we also write as
$$
\partial _s^lP^{-s}={\Cal P}_l(P)P^{-s}, \text{ with } {\Cal P}_l(P)=(-\log P)^{l}.\tag3.10
$$

The value of $\partial _s^l\zeta (I,P,s)$ at $s=0$, more generally
the constant term of 
\linebreak $\partial _s^l\zeta (A,P,s)$ at $s=0$, will be
determined as
a specific coefficient in expansions of  
$
\Gamma (s)\Tr( A{\Cal P}_l(P) P^{-s})$ and $\Tr ( A{\Cal P}_l(P) (P-\lambda )^{-N})$.
A framework for such calculations has been set up in Lesch \cite{L}.
With the
notation for symbol spaces introduced there, $A{\Cal P}_l(P)$ is log-polyhomogeneous
belonging to $\operatorname{CL}^{\nu 
,l}(X)$. 
In local coordinates, the symbols of such operators have the
structure$$
b(x,\xi )\sim\sum_{j\in\Bbb N}\sum_{\sigma =0}^lb_{\nu
-j,\sigma }(x,\xi )\log^\sigma [\xi ],\tag3.11
$$
with $b_{\nu -j,\sigma }$ homogeneous in $\xi $ of degree $\nu -j$
for $|\xi |\ge 1$.
This defines the symbol space $\operatorname{CS}^{\nu 
,l}$. The symbols  (and the operators they define) are said to be of
order $\nu $; 
the degree of a term is the number $\nu -j$.
Log-polyhomogeneous operators were studied earlier by Schrohe \cite{Sc}.

A generalization of \cite{GS1, Th\. 2.1} to such operators is shown
in \cite{L, Th\. 3.7}, which we develop further in Theorem 3.2
below.
For this, we recall from \cite{L} that the definition of the finite
part integral $\tslint
f(x,\xi )\,\d\xi $ in (1.18), (1.35) can be extended to the symbols
$f(x,\xi )\in 
\operatorname{CS}^{\nu ,l}(\Bbb R^n)$:  Omit $x$-dependence.
When $\nu <-n$, $\tslint
f(\xi )\,\d\xi $ is the usual integral $\int_{\Bbb R^n}f(\xi
)\,\d\xi $; more generally
it is equal to the constant term
$p_0(0)$ in the
asymptotic expansion of $\int_{|\xi |\le \mu }f(\xi )\,\d\xi $ in
powers and log-powers of $\mu $:
$$
\aligned
\slint
f(\xi )\,\d\xi &=p_0(0), \text{ when}\\ 
\int_{|\xi |\le \mu }f(\xi )\,\d\xi &\sim \sum_{j\in\Bbb N, j\ne \nu
+n}p_{\nu +n-j}(\log \mu )\mu ^{\nu +n-j}+p_0(\log \mu )\mu ^0,\text{ for
}\mu \to\infty .
\endaligned\tag3.12
$$
Here the $
p_{\nu +n-j}$ are polynomials of degree $\le l$ when $\nu +n-j\ne 0$,
degree $\le l+1$ when   $\nu +n-j= 0$. 
An explicit formula
for $p_0(0)$ defined from a log-homogeneous term is worked out
in \cite{L, (5.12)}:  When $f(x,\xi )=f_\nu (x,\xi )\log^\sigma [\xi ]$
with $f_\nu $  homogeneous of degree $\nu $ in $\xi $ for $|\xi |\ge
1$, then$$ 
\slint f(x,\xi )\,\d\xi =\int_{|\xi |\le 1}f(x,\xi )\,\d\xi +\tfrac{(1-\delta
_{\nu +n,0})
(-1)^{\sigma +1}\sigma !}{(\nu +n)^{\sigma +1}}\int_{|\xi
|=1}f_{\nu }(x,\xi )\,\d S(\xi ).
\tag3.13
$$

It follows that when $b$ is as in (3.11), then 
$$
\multline
\slint b(x,\xi )\d \xi \\
=\sum_{j\le \nu +n}\sum_{0\le\sigma \le l}
\bigl(\int_{|\xi |\le 1} b_{\nu -j,\sigma }(x,\xi )\,\d\xi +\tfrac{(1-\delta
_{\nu +n, j})
(-1)^{\sigma +1}\sigma !}{(\nu +n-j)^{\sigma +1}}\int_{|\xi
|=1}b_{\nu -j,\sigma }(x,\xi )\,\d S(\xi )
\bigr)\\
+\int_{\Bbb
R^{n}}
\bigl(b(x,\xi )-\sum_{j \le \nu +n}\sum_{0\le\sigma \le l}b_{\nu -j,\sigma
}(x,\xi )\log^\sigma [\xi ]\bigr)\,\d\xi. 
\endmultline\tag 3.14
$$
As in (1.18), the sum over $j\le \nu +n$ can be replaced by the sum
over $j\le J$ for any choice of $J\ge \nu +n$.

The definition of having {\it even-even resp.\ even-odd
alternating 
parity} is extended to symbols (3.11) to mean that
$$
\aligned
\text{even-even: }b_{\nu -l,\sigma }(x,-\xi )&=(-1)^{\nu -l}b_{\nu -l,\sigma }(x,\xi )
\text{ for }|\xi |\ge 1,\text{ resp.}\\
\text{even-odd: }b_{\nu -l,\sigma }(x,-\xi )&=(-1)^{\nu -l-1}b_{\nu
-l,\sigma }(x,\xi ) \text{ for }|\xi |\ge 1,
\endaligned \tag3.15 
$$
with similar properties of the derivatives.
Then we can consider the four cases in Definition 1.1 for
log-polyhomogeneous operators $B$. 

\proclaim{Theorem 3.2} {\rm (i)} Let $\nu \in\Bbb R$ and $l\in \Bbb N$, let
$B$ be log-polyhomogeneous in 
$\operatorname{CL}^{\nu 
,l}$ with symbol {\rm (3.11)} on $\Bbb R^n$. Let $P$ be a classical
$\psi $do, uniformly elliptic of integer order $m>0$ on $\Bbb R^n$ and
with no principal symbol eigenvalues in a sector around $\Bbb R_-$.
Let $N>(\nu +n)/m$. There is an
asymptotic expansion of the kernel of $B(P-\lambda )^{-N}$ at $x=y$:$$
K( B (P-\lambda )^{-N},x,x)\sim\sum_{j=0}^\infty
\sum_{\sigma=0}^{l+1}\tilde c_{j,\sigma}(x)(-\lambda )^{\frac{\nu
+n-j}m-N}\log^{\sigma}(-\lambda ) +\sum_{k=0}^\infty 
\tilde c''_{k}(x)(-\lambda )^{-k-N};\tag3.16
$$
here $\tilde c_{j,l+1}(x)=0$ unless $\frac{j-\nu -n}m\in \Bbb N$. 
The coefficients $\tilde c_{j,\sigma}(x)$
depend on the homogeneous or log-homogeneous symbols of the
first $j+1$ degrees in $B$ and $P$ (are local in this sense). The
$\tilde c^{\prime\prime}_k$ depend on the full structure (are global).
In particular, when we define $\tilde c_{\nu +n,0}(x)=0$ if $\nu <-n$ or $\nu \in \Bbb
R\setminus \Bbb Z$,
$$
\tilde c_{\nu +n,0}(x)+\tilde c''_0(x)=\slint b(x,\xi )\,\d\xi
+\text{ local terms}.\tag 3.17
$$

{\rm (ii)} It follows that in the comparison of the coefficients
for two choices of auxiliary operator $P$ and $P'$, $
\tilde c''_{0}(B,P,x)-\tilde c''_0(B,P',x)$ is local (in the above sense).

{\rm (iii)} If, moreover,  $B$ and $P$ satisfy {\rm (3)} or {\rm (4)}
of Definition {\rm 1.1} (in particular,  $m$ is even), then the expansion {\rm (3.16)} reduces
to 
the form$$\multline
K(B(P-\lambda )^{-N},x,x)\\
\sim\sum_{j\in\Bbb N, j-n-\nu \,\text{\rm odd\;}}
\sum_{\sigma=0}^{l}\tilde c_{j,\sigma}(x)(-\lambda )^{\frac{\nu
+n-j}m-N}\log^{\sigma}(-\lambda ) +\sum_{k=0}^\infty 
\tilde c''_{k}(x)(-\lambda )^{-k-N};\endmultline\tag3.18
$$

{\rm (iv)} In each of the cases {\rm (1)--(4)} of Definition {\rm
1.1}, the $\tilde c_{\nu +n,\sigma }(x)$ vanish and $$
\tilde c_0''(x)=\slint b(x,\xi )\,\d\xi .\tag 3.19
$$

{\rm (v)} If the symbol of $B$ is given in $(x,y)$-form $b(x,y,\xi
)$, the formulas hold with 
$b(x,\xi )$ replaced by $b(x,x,\xi )$.
\endproclaim  

\demo{Proof} The proof of \cite{L, Th\. 3.7} is in fact modeled
very closely after the proof of \cite{GS1, 
Th\. 2.1}, which we recalled to a large extent in the proof of Theorem 1.3
above. Let $$ 
Q(\mu )=B(P-\lambda )^{-N},\text{ with }-\lambda =\mu ^m.\tag3.20
$$
The symbol $q(x,\xi ,\mu )$ is now the composite of a
log-polyhomogeneous symbol  in 
$\operatorname{CS}^{\nu ,l}$ and a weakly polyhomogeneous symbol in
$S^{-mN,0}\cap S^{0,-mN}$. Here, since $\operatorname{CS}^{\nu
,l}\subset S^{\nu +\varepsilon }$ for any $\varepsilon >0$, $$
q(x,\xi ,\mu )\in S^{\nu +\varepsilon -mN,0}\cap S^{\nu +\varepsilon
,-mN}
\tag3.21$$
(defined also for $S_{1,0}$-symbols without requirements of
polyhomogeneity). Now 
there are expansions   $$
q(x,\xi ,\mu )=\mu ^{-mN}\sum_{0\le l'<L}q^{(l')}(x,\xi )\mu ^{-ml'}+O(\ang\xi
^{\nu +\varepsilon +mL}\mu ^{-m(N+L)}),\tag3.22
$$
for all $L$, with coefficient symbols $q^{(l')}(x,\xi )$ 
log-polyhomogeneous in 
$\operatorname{CS}^{\nu +ml',l}$. One analyzes the kernel defined from
$q$ by (1.24), by splitting the contribution from each log-homogeneous
term into three pieces as in (1.25). The integral over $\{|\xi |\ge
\mu \}$ contributes to the $\tilde c_{j,\sigma}$-terms. For the
integral over $\{|\xi |\le 1\}$ one uses (3.22) and gets contributions
to the $\tilde c^{\prime\prime}_{k}$-terms. For the integral over $\{1\le |\xi
|\le \mu \}$ one likewise uses the expansion (3.22) for each
log-homogeneous 
term in the symbol; each expansion term gives a contribution$$
\aligned
&\mu ^{-m(N+l')}\int_{1\le |\xi |\le \mu }q^{(l')}_{\nu -j+ml',\sigma }(x,\xi 
)\log^\sigma |\xi |\,\d\xi \\
&=\mu ^{-m(N+l')}\int_1^\mu r^{\nu -j+ml' +n-1}\log^\sigma r\,dr
\int_{|\xi |=1}q_{\nu -j+ml',\sigma }^{(l')} (x,\xi)\,\d S(\xi )\\
&=\cases \sum_{\sigma '=0}^\sigma c_{j,l',\sigma '}\mu ^{\nu -j-mN
+n}\log^{\sigma '}
\mu +c''_{j,l'}\mu ^{-m(N+l')}&\text{ if 
}\nu -j+ml'+n\ne 0,\\
c'_{j,l'}\mu ^{-m(N+l')}\log ^{\sigma +1}\mu &\text{ if }\nu
-j+ml'+n= 0.\endcases 
\endaligned\tag3.23
$$
The coefficients $c_{j,l',\sigma '}$ and $c'_{j,l'}$ contribute to the 
$\tilde c_{j,\cdot}$-terms, whereas the coefficient $c''_{j,l'}$
contributes to the $\tilde c^{\prime\prime}_{l'}$-term. They are
proportional to $\int_{|\xi |=1}q_{\nu -j+ml',\sigma }^{(l')} \,\d S$
by universal factors; the value of $c''_{j,l'}$ is 
$$
c''_{j,l'}(x)=\frac{(-1)^{\sigma +1}\sigma !}{(\nu -j+ml'+n)^{\sigma
+1}}
\int_{|\xi |=1}q_{\nu -j+ml',\sigma }^{(l')} (x,\xi)\,\d S(\xi ),\tag3.24
$$
cf\. \cite{L, (5.12)}.  (This term was left
out in \cite{L, (3.38)}; the connection between $\tilde c''_0$
and $\tslint b$
was made only towards the end of the paper.) 
Remainders are treated essentially as in
\cite{GS1}, as recalled above in Theorem 1.3; this shows (3.16). 

In the cases (3) and (4) as in Definition 1.1, the terms in the sum
over $j$ vanish for $j-n-\nu $ even, since they are obtained by integration
in $\xi $ of odd functions (like in Theorem 1.3); this shows
(iii).

The analysis leading to (3.17) is
practically the same as in the proof of (1.17) in Theorem 1.3, only
with $a$ replaced by $b$; $P$ is unchanged. 
It is seen again that
all parts of $(P+\mu ^m)^{-N}$ except $\operatorname{OP}((p_m(y,\mu
)+\mu ^m)^{-N})$ gives series with a 
locally determined constant term, (1.30) being replaced by$$
K( x,x,\mu )\sim\sum_{j=0}^\infty
\sum_{\sigma=0}^{l+1}\tilde d_{j,\sigma}(x)(-\lambda )^{\frac{\nu
+n-j}m-N}\log^{\sigma}(-\lambda ) +\sum_{k=1}^\infty 
\tilde d''_{k}(x)(-\lambda )^{-k-N}.\tag3.25
$$
The considerations on the integral $\int_{|\xi |\le \mu }b(x,\xi
)(p_m(x,\xi )+\mu ^m)^{-N}\,\d\xi $ carry over vebatim from the
considerations on $\int_{|\xi |\le \mu }a(x,\xi
)(p_m(x,\xi )+\mu ^m)^{-N}\,\d\xi $ in Theorem 1.3. 

This shows
(3.17), and (ii) is an immediate consequence since the $\tilde
c_{\nu +n,0}$ are local and the symbol integrals cancel out.
Moreover, (iv) is seen by observing that the local 
contributions (from integrals over $|\xi |\ge \mu $ and from the sum
over $j$ in the various series of the 
form (3.25) that arise in the analysis) vanish in the cases (1)--(4) of
Definition 1.1. 

Finally, (v) is included as in Remark 1.4.
\qed
\enddemo 

It is known from \cite{L} in the cases (1) and (2) for
log-polyhomogeneous operators that the density $\tslint b(x,\xi
)\,\d\xi \,|dx|$ or $\tslint b(x,x,\xi )\,\d\xi \,|dx|$ has an
invariant meaning, 
and the argument
carries
over to log-polyhomogeneous operators in the parity cases (3) and
(4), in the same way as mentioned after Theorem 1.3.
So, in the cases (1)--(4), when $B $ is given on $X$, it defines a density
$\omega_{\operatorname{TR}}(B)$
described in local coordinates by $$
\omega_{\operatorname{TR}}(B)=\slint B(x,\xi )\,\d\xi |dx|\;\text{
resp\. }\;\omega _{\operatorname{TR}}(B)=\slint b(x,x,\xi )\,\d\xi\,
|dx|,\tag 3.26  
$$
when $B$ has symbol in $x$-form $b(x,\xi )$, resp\. in  $(x,y)$-form
$b(x,y,\xi )$. 

The inclusion of symbols in $(x,y)$-form allows us in
particular to observe that when $B_1B_2$ is as in one of the cases (1)--(4), 
with $B_1=\operatorname{OP}(b_1(x,\xi
))$ and $B_2=\operatorname{OP}(b_2(y,\xi ))$ in a local coordinate
system, then$$
\omega _{\operatorname{TR}}(B_1B_2)=\slint b_1(x,\xi )b_2(x,\xi
)\,\d\xi\,|dx| .\tag3.27$$ 

We have as usual a corollary on the manifold situation, when $B$ is
decomposed as in (1.13) 
and the pieces are carried over to local coordinates in $\Bbb R^n$ as
explained there.

\proclaim{Corollary 3.3} Consider a log-polyhomogeneous operator $B$
on the manifold $X$, together with $P$ as in Proposition {\rm 2.1}, 
with $N>(\nu +n)/m$. Then there is an asymptotic
expansion of the trace:
$$
\Tr( B (P-\lambda )^{-N})\sim\sum_{j\in\Bbb N}
\sum_{\sigma=0}^{l+1}\tilde c_{j,\sigma}(-\lambda )^{\frac{\nu
+n-j}m-N}\log^{\sigma}(-\lambda ) +\sum_{k=0}^\infty 
\tilde c''_{k}(-\lambda )^{-k-N}.\tag3.28
$$
Here,
when we define $\tilde c_{\nu +n,0}=0$ if $\nu <-n$ or $\nu \in \Bbb
R\setminus \Bbb Z$, and consider the operators localized to $\Bbb R^n$
as explained before Theorem {\rm 1.3}, then
$$
\tilde c_{\nu +n,0}+\tilde c''_0=\int_{\Bbb
R^n}\slint \tr b(x,\xi )\,\d\xi dx +\text{ local terms};\tag3.29
$$
where $ b(x,\xi )$ is the symbol of $B$; 
the local terms depend only on the 
strictly homogeneous
terms in the symbols of $B$ and $P$ for $j\le n+[\nu]$. (If the
symbol $b$ is
in $(x,y)$-form, the formula holds with $b(x,x,\xi )$ instead.)  

{\rm (ii)} It follows that in the comparison of the coefficients
for two choices of auxiliary operator $P$ and $P'$, $
\tilde c_{\nu +n,0}(B,P)+\tilde c''_{0}(B,P)-(\tilde c_{\nu +n,0}(B,P')+\tilde c''_0(B,P'))$ is local.

{\rm (iii)} If, moreover,  $B$ and $P$ satisfy
{\rm (3)} or {\rm (4)} of Definition {\rm 1.1}, then the expansion reduces
to 
the form
$$
\Tr( B (P-\lambda )^{-N})\sim\sum_{j\in\Bbb N, j-n-\nu \,\text{\rm odd\;}}
\sum_{\sigma=0}^{l}\tilde c_{j,\sigma}(-\lambda )^{\frac{\nu
+n-j}m-N}\log^{\sigma}(-\lambda ) +\sum_{k=0}^\infty 
\tilde c''_{k}(-\lambda )^{-k-N}.\tag3.30
$$
In particular, $\zeta (B,P,s)$ (cf.\ {\rm (3.32)} below) is regular
at all integers $s$.

{\rm (iv)} In each of the cases {\rm (1)--(4)} in Definition {\rm
1.1}, the $\tilde c_{\nu +n,\sigma }$ vanish (for any choice of local
coordinates), and (cf\. {\rm (3.26)})$$
\tilde c_0''=\int _X\tr \omega _{\operatorname{TR}}(B).\tag 3.31
$$
\endproclaim 

In the following, we draw on the hypothesis that $P$ is invertible.
By the transition formulas in \cite{GS2}, (3.28) implies the structure of
the meromorphic
extension of $\Gamma (s)\Tr( B P^{-s})$, also denoted $\Gamma
(s)\zeta (B,P,s)$: 
$$
\Gamma (s)\Tr( B P^{-s})=\Gamma
(s)\zeta (B,P,s)\sim\sum_{j=0}^\infty
\sum_{\sigma=0}^{l+1}\frac{ c_{j,\sigma}}{(s+\frac{j-\nu
-n}m)^{\sigma+1}} 
+\sum_{k=0}^\infty \frac {c^{\prime\prime}_{k}}{s+k}
,\tag3.32
$$
with universal nonzero factors linking $\tilde c_{ j,\sigma}$ with 
$c_{ j,\sigma}$ and $\tilde c^{\prime\prime}_{ k}$ with 
$c^{\prime\prime}_{ k}$; in particular, $c''_0=\tilde c''_0$ and
$c_{\nu +n,0}=\tilde c_{\nu +n,0}$ (with the usual zero convention if
the series in $j$ does not contain such a term). Dividing out $\Gamma
(s)$, we see that  
$\zeta (B,P,s)$ has a Laurent expansion at $s=0$ when $B\in
\operatorname{CL}^{\nu ,l}(X)$:  
$$
\zeta (B,P,s)\sim \sum_{r\ge -l-1}C_{r}(B,P)s^{r};\tag 3.33
$$ 
here $$
C_0(B,P)=c_{\nu +n,0}+c''_0=\tilde c_{\nu +n,0}+\tilde c''_0.\tag3.34
$$ 

Corollary 3.3 (ii) shows that $C_0(B,P)$ satisfies
the first 
condition in Definition 2.2 for being a quasi-trace.
(In \cite{O1}, the difference $C_0(B,P)-C_0(B,P')$ is shown to be
a certain residue, in the case where the symbol of $B$ equals
$c\log[\xi ]$ plus a zero-order classical symbol.)  We shall now 
show that the second condition, concerning commutators, is likewise
satisfied.

\proclaim{Theorem 3.4} 
Let $B\in \operatorname{CL}^{\nu ,l}$, $B'\in
\operatorname{CL}^{\nu ',l'}$,  
and
let $P$ be as in Proposition
{\rm 2.1 (ii)}. Then
$$
\Tr( [B ,B'] (P-\lambda )^{-N})\sim\sum_{j=0}^\infty
\sum_{\sigma=0}^{l+1}\tilde b_{j,\sigma}(-\lambda )^{\frac{\nu +\nu '
+n-j}m-N}\log^{\sigma}(-\lambda ) +\sum_{k=1}^\infty \tilde
b''_{k}(-\lambda )^{-k-N}.\tag3.35
$$ Hence the zeta function $\zeta ([B ,B'],P,s)=\Tr( [B ,B'] P^{-s})$
satisfies (when $P$ is invertible)
$$
\zeta ([B,B'],P,s)\sim \sum_{r\ge -l+l'-1}C_{r}([B,B'],P)s^{r}\text{
for }s\to 0,\tag 3.36
$$ 
with $C_0([B,B'],P)$ local, 
depending solely on the terms of the first 
$n+[\nu +\nu ']+1$
homogeneity degrees in the symbols of $B $, $B'$ and $P$.

In particular, it vanishes if one of the conditions {\rm
(1$'$)--(4$'$)} in Theorem {\rm
1.7} is satisfied (with $A$, $A'$ replaced by $B$, $B'$).
\endproclaim 

\demo{Proof} Since $[B,B']$ has symbol in $\operatorname{CS}^{\nu
+\nu ',l+l'}$, we know already from Corollary 3.3 that
$Q=[B,B'](P-\lambda )^{-N}$ has the corresponding  expansion 
(3.28); we just have to show that $\tilde c''_0=0$ there.

As in (1.47) and (1.48), the trace calculation can be
reduced to the calculation for
$$
Q(\mu )=B \sum_{0\le M<N}c_{MN}(P+\mu ^m)^{-1-M}[P,B'](P+\mu
^m)^{-N+M};\tag3.37
$$
by the rules of calculus, it has symbol in $S^{\nu
+\nu '-mN+\varepsilon ,0}\cap S^{\nu +\nu '+m+\varepsilon ,
-(m+1)N}$. The trace expansion can be analyzed as in the proof of
Theorem 3.2;
with the modification
that the contributions from the integral over $\{|\xi |\le 1\}$ start
with the lower power $\mu ^{-m(N+1)}$, and the terms of the form
$-c''_{j,l'}\mu ^{-\alpha }$ from (3.23) all have $\alpha \ge m(N+1)$.
Thus $\tilde c''_0=0$ for this operator family. In the special cases
(1$'$)--(4$'$), also the contribution from the series in $j$
vanishes, as shown in Corollary 3.3.
\qed
\enddemo 

We have hereby obtained:

\proclaim{Corollary 3.5} 
{\rm (i)} For log-polyhomogeneous operators $B$ with $P$ as in
Proposition {\rm 2.1}, $C_0(B,P)$ {\rm (3.34)} is a quasi-trace in the
sense of Definition {\rm 2.3}.

{\rm (ii)} Definition {\rm 1.1} of the canonical trace
$\operatorname{TR}$ extends to
log-polyhomogeneous operators $B$, in such a way that the 
definition is independent of
$P$, Theorem {\rm 1.7} extends to these operators, and$$
\operatorname{TR}(B)=\int_{X}\tr \omega_{\operatorname{TR}}(B).\tag3.38
$$
\endproclaim 

Point (ii) was shown in \cite{L, Sect\. 5} for the cases where 
$\nu \notin\Bbb Z$ or $<-n$, by a somewhat different proof
(where $P$ was assumed to be selfadjoint positive with scalar principal
symbol). 
\medskip

Let us now return to the zeta function for a classical $\psi $do,
cf\. (3.1). The expansions
(3.28) and (3.32)
hold in particular for $B=A{\Cal P}_l(P)$; let us denote the coefficients 
in this case by $\tilde c^{(l)}_{j,\sigma}$,  $\tilde
c^{\prime\prime(l)}_{k}$,
resp\. $ c^{(l)}_{j,\sigma}$,  $ c^{\prime\prime(l)}_{k}$.
Then we have found that$$
\partial _s^l\zeta (A,P,s)=\Tr(A{\Cal P}_l(P)P^{-s})=\zeta (A{\Cal P}_l(P),P,s)\tag3.39
$$
has
the meromorphic structure determined from 
$$
\Gamma (s)\Tr( A{\Cal P}_l(P) P^{-s})\sim\sum_{j=0}^\infty
\sum_{\sigma=0}^{l+1}\frac{ c^{(l)}_{j,\sigma}}{(s+\frac{j-\nu
-n}m)^{\sigma+1}} 
+\sum_{k=0}^\infty \frac {c^{\prime\prime(l)}_{k}}{s+k}
.\tag3.40
$$

Concerning $\zeta (A,P,s)$ in (3.1), we find by
differentiation in $s$:
$$
\zeta (A{\Cal P}_l(P),P,s)=\partial _s^l\zeta
(A,P,s)\sim\tfrac{(-1)^ll!}{s^{l+1}}C_{-1}(A,P)+\sum_{r\ge l}\tfrac{r!}{(r-l)!}C_r(A,P)s^{r-l},\tag3.41
$$
so that $C_0(A{\Cal P}_l(P),P)=l!C_l(A,P)$ for $l\ge 1$.
Dividing by $\Gamma (s)$ in (3.40), we see that the constant term 
at $s=0$ in the expansion of $\zeta ( A{\Cal P}_l(P),P,s)$  is
$
c^{(l)}_{\nu +n,0}+c^{\prime\prime(l)}_{0}$, with the usual
conventions.
So we have found:

\proclaim{Corollary 3.6}
For $l\ge 1$, $$
C_l(A,P)=\tfrac1{l!}C_0(A{\Cal P}_l(P),P)=\tfrac 1{l!}(c^{(l)}_{\nu +n,0}+c^{\prime\prime(l)}_{0}),\tag3.42
$$
defined from the constants appearing in {\rm (3.40)}; here
$c^{(l)}_{\nu +n,0}=0$ if $\nu +n\notin \Bbb N$ or $\nu <-n$.
\endproclaim

Note that $c^{(l)}_{\nu +n,0}$ is locally determined (in local
coordinates it 
comes from the part of
the symbol of $A{\Cal P}_l(P)(P+\mu ^m )^{-N}$ with homogeneity degree $
-mN-n$), whereas $c^{\prime\prime(l)}_{0}$ depends on the full structure (is
global). 

It may also be observed that the coefficient of $s^{-l-1}$ in (3.41),
proportional to $\operatorname{res} A$, is also
proportional to the $(l+1)$'st higher residue
of $A{\Cal P}_l(P)$ defined in \cite{L, Sect\. 4}.

In particular, $
C_1(A,P)$ equals $c^{(1)}_{\nu +n,0}+c^{\prime\prime(1)}_{0}$;
this specializes to a formula for $-\log\det P$ in the case $A=I$:

\proclaim{Corollary 3.7} One has that$$
-\log\det P =C_1(I,P)=C_0(-\log P ,P)=c^{(1)}_{\nu
+n,0}+c^{\prime\prime(1)}_{0},\tag3.43 
$$
determined from {\rm (3.40)} in the case $A=I$.
\endproclaim

One cannot conclude from Corollary 3.6 that the higher Laurent
coefficients in (3.1) are quasi-traces of $A$ itself --- for in $C_0(A{\Cal P}_l(P),P)$,
the first entry depends highly on the choice of $P$ when $l>0$.
However, when parity and dimension match, these coefficients can
be expressed by the extended $\operatorname{TR}$ applied to
$A{\Cal P}_l(P)$. 
In fact, when $P$ is even-even, ${\Cal P}_l(P)$ is even-even for all $l\ge 1$.
Then one gets, by Corollary 3.6:

\proclaim{Corollary 3.8} Let $P$ be even-even. If $n$ is odd, let the
classical $\psi $do $A$
be even-even; if $n$ is even, let $A$ be even-odd. Then$$
C_l(A,P)=\tfrac1{l!}C_0(A{\Cal P}_l(P),P)=\tfrac1{l!}\operatorname{TR}(A{\Cal P}_l(P)),
\text{ for }l\ge 1.\tag3.44
$$
In particular:

{\rm (i) } When $n$ is odd, then 
$$\aligned
\log\det P&=\operatorname{TR}(\log P),\\
\partial _s^l\zeta (I,P,0)&=\operatorname{TR}({\Cal P}_l(P)),\text{ for $l\ge 1$}.
\endaligned\tag3.45$$

{\rm (ii) } When $n$ is even and $D$ is a first-order
selfajoint elliptic differential operator (cf\. Remark {\rm 1.4}), 
$$
\partial _s^l\eta
(D,0)=2^{-l}\operatorname{TR}(D|D|^{-1}{\Cal P}_l(D^2)),\text{ for $l\ge 1$}. 
\tag3.46$$

\endproclaim 

A formula similar to the first line in (3.45) appears in the abstract
of \cite{O3} (it can 
be justified on the basis of \cite{O2}, Lemma 0.1).
As an example of the second line, 
$$
\partial _s^2\zeta (I,P,0)=\operatorname{TR}((\log P)^2).\tag 3.47
$$
cf\. (3.10). See also \cite{KV, Sect\. 4}.

\example{Remark 3.9} Let us set the above methods in relation to the
results of \cite{O1}, \cite{O2} and \cite{KV} on the
multiplicative anomaly $\det AB/(\det A \det B)$ of the determinant. 
\cite{O1} shows 
that when the elliptic positive-order $\psi $do's $A$, $B$ and $AB$ have
scalar principal symbol taking no values on $\rmi$, and have no
eigenvalues on $\crm$, then$$
\log AB -\log A-\log B=[A, C^{(I)}(A,B)]+[B, C^{(II)}(A,B)] +F,\tag3.48
$$where $C^{(I)}(A,B)$   and $C^{(II)}(A,B)$ are Lie polynomials in
$A$ and $B$, and $F$ is an operator of order $<-n$ with $\Tr F=0$. In
\cite{O2} this is used to show that 
$\log \det AB-\log\det A-\log\det B$ 
equals the noncommutative residue of a certain operator derived
from $A$ and $B$; in particular it is locally determined. By
variational methods, \cite{O2} and \cite{KV} show local determinedness of 
$\log \det AB-\log\det A-\log\det B$ also in cases where the
principal symbol is not scalar. This implies local
determinedness of $\det AB/(\det A \det B)$ (by exponentiation). 

From our point of view, (3.48) implies that $C_0(\log 
AB-\log A-\log B, P)$ is locally determined for any auxiliary
operator $P$, by Theorem 3.4 and the fact that $\Tr F=0$. Then
since the expressions $C_0(\log A,A)-C_0(\log A,P)$, $C_0(\log B,B)-C_0(\log B,P)$
and $C_0(\log AB,AB)-C_0(\log AB,P)$ are locally determined by
Corollary 3.3 (ii), we conclude the local determinedness of $$
\log \det AB-\log\det A-\log\det B=C_0(\log AB,AB) -C_0(\log A,A)
-C_0(\log B,B).  
$$

When $n$ is odd and $A$ and $B$ are even-even, then (cf\.
Corollary 3.6 (i)) 
$$\log \det AB-\log\det A-\log\det B=\operatorname{TR}(\log AB-\log
A-\log B)=0,$$ since it is locally determined (local contributions
give zero because of parity), so $\det AB=\det 
A\det B$ then, as originally shown in \cite{KV, Th\. 7.1}.
\endexample

\enddocument